\documentclass[12pt]{amsart}

\usepackage{geometry}
\newcommand{\ye}[1]{{\textcolor{purple}{[Ye: #1]}}}

\usepackage{amsfonts, adjustbox, amssymb, amscd}
\numberwithin{equation}{section}

\usepackage{bm}
\usepackage{verbatim}
\usepackage{mathrsfs}
\usepackage{graphicx}
\usepackage{tikz-cd}
\usepackage{subcaption}
\usepackage{listings}
\usepackage{subfiles}
\usepackage[toc,page]{appendix}
\usepackage{mathtools}
\usepackage{comment}
\usepackage{enumerate}
\usepackage{enumitem}
\usepackage[all]{xy}

\usepackage{graphicx}
\usepackage{appendix}
\usepackage{hyperref}
\hypersetup{
    colorlinks=true,
    citecolor=red,
    linkcolor=blue,
    filecolor=magenta,      
    urlcolor=red,
}
\newcommand{\bir}{\dashrightarrow}

\newcommand{\bb}{\bm{b}}
\newcommand{\Mm}{{\bf{M}}}
\newcommand{\Nn}{{\bf{N}}}

\newcommand{\Pp}{{\bf{P}}}

\newcommand{\Dd}{{\bf{D}}}

\newcommand{\Spec}{\mathrm{Spec}}
\newcommand{\id}{\mathrm{id}}

\newcommand{\Qq}{\mathbb{Q}}

\newcommand{\Rr}{\mathbb{R}}

\newcommand{\Zz}{\mathbb{Z}}

\newcommand{\ww}{\mathbf{w}}
\newcommand{\vv}{\bm{v}}
\newcommand{\uu}{\mathbf{u}}

\newcommand{\Kk}{\mathbb{K}}

\newcommand{\Center}{\operatorname{center}}

\newcommand{\Exc}{\operatorname{Exc}}

\newcommand{\Conv}{\operatorname{Conv}}

\newcommand{\Nklt}{\operatorname{Nklt}}
\newcommand{\mld}{{\operatorname{mld}}}

\newcommand{\lct}{\operatorname{lct}}

\newcommand{\Supp}{\operatorname{Supp}}

\newcommand{\mult}{\operatorname{mult}}

\newcommand{\Bb}{{\bf{B}}}

\newcommand{\Ii}{\Gamma}

\newcommand{\Ss}{\mathcal{S}}

\newcommand{\Nef}{\mathrm{Nef}}


\newcounter{parentnumber}

\newtheorem{thm}{Theorem}[section]
\newtheorem{conj}[thm]{Conjecture}
\newtheorem{cor}[thm]{Corollary}
\newtheorem{lem}[thm]{Lemma}
\newtheorem{prop}[thm]{Proposition}

\theoremstyle{definition}
\newtheorem{defn}[thm]{Definition}

\theoremstyle{definition}
\newtheorem{rem}[thm]{Remark}

\newtheorem{deflem}[thm]{Definition-Lemma}

\newtheorem{nota}[thm]{Notation}

\theoremstyle{definition}

\begin{document}

\title[Existence of minimal models for 3-fold g-pairs in char $p>0$]{Existence of minimal models for threefold generalized pairs in positive characteristic}
\author{Tianle Yang}
\address{School of Mathematical Sciences, East China Normal University, 500 Dongchuan Road, 200241 Shanghai, China}
\email{10211510098@stu.ecnu.edu.cn}

\author{Zelin Ye}
\address{Shanghai Center for Mathematical Sciences, Fudan University, 2005 Songhu Road, 200438 Shanghai, China}
\email{zlye24@m.fudan.edu.cn}

\author{Zhiyao Zhang}
\address{Department of Mathematics, Shanghai Jiao Tong University, 800 Dongchuan Road, 200240 Shanghai, China}
\email{aboctopus@sjtu.edu.cn}

\subjclass[2020]{14E30, 14B05.}
\keywords{Minimal model program. Positive characteristic. Threefolds.}
\date{\today}

\begin{abstract}
Let $\Kk$ be an algebraically closed field of characteristic $p>5$. We show the existence of minimal models for pseudo-effective NQC lc generalized pairs in dimension three over $\Kk$. As a consequence, we prove the termination of flips for pseudo-effective threefold NQC lc generalized pairs over $\Kk$. This provides a new proof of the termination of flips for pseudo-effective pairs over $\Kk$ without using non-vanishing theorems. A key ingredient of our proof is the ACC for lc thresholds in dimension $\leq 3$ and the global ACC in dimension $\leq 2$ for generalized pairs over $\Kk$.
\end{abstract}

\maketitle

\tableofcontents

\section{Introduction}\label{sec:Introduction}

Throughout the paper, we work, without further notice, over an algebraically closed field $\mathbb K$ of characteristic $p > 5$.

Over any algebraically closed field of characteristic $0$, \cite{BCHM10} proved the existence of klt flips, thereby establishing the existence of the minimal model program over any algebraically closed field of characteristic $0$. Following this groundbreaking work, there have been significant developments in minimal model program theory over the past two decades. In particular, much work has been dedicated to improving or varying the corresponding results in \cite{BCHM10}. These improvements and variations usually focus on three main directions: relaxing the restriction on singularities (e.g., from klt singularities to lc singularities, cf. \cite{Bir12, HX13, HH20}), considering new structures (e.g., generalized pairs, cf. \cite{BZ16, HL23}, foliations, cf. \cite{CS21, CS25, CHLX23}), and changing the base field (e.g., positive characteristic, cf. \cite{Bir16, HX14, HNT20, HW22}, mixed characteristic, cf. \cite{BMPSTWW23}).

It is interesting to ask whether these generalizations can be combined, e.g., by considering multi-directional rather than one-directional variations. For example, \cite{HL23} considered the minimal model program for generalized pairs with lc singularities, while \cite{HNT20} considered the minimal model program for lc pairs over fields of characteristic $p>5$. In this paper, we consider the following variation: the minimal model program for (klt) generalized pairs over an algebraically closed field of positive characteristic, or, more precisely, of characteristic $p>5$.

There has already been some progress on the minimal model program for klt generalized pairs over fields of positive characteristic \cite{BF25, BS23, FW23}. In particular, the existence of the minimal model program is known in dimension $\leq 3$. However, many fundamental theorems of the minimal model program for generalized pairs over fields of positive characteristic remain unknown in dimension $\leq 3$, particularly the existence of minimal models.

The goal of this paper is to systematically study the minimal model program for (NQC) klt generalized pairs in dimension $3$ over an algebraically closed field of characteristic $p > 5$. The first main theorem of our paper is the following:

\begin{thm}[Termination of pseudo-effective flips]\label{thm: tof intro}
Let $(X, B,\Mm)/U$ be an NQC lc generalized pair of dimension $\leq 3$ over an algebraically closed field of characteristic $p>5$ such that $K_X+B+\Mm_X$ is pseudo-effective$/U$. Then any sequence of $(K_X+B+\Mm_X)$-flips$/U$ terminates.
\end{thm}

Note that when $\Mm = \bm{0}$, Theorem \ref{thm: tof intro} becomes the termination of flips for pseudo-effective lc pairs, which was proven in \cite{Xu24} and essentially relies on the non-vanishing theorem for lc pairs in dimension $3$. Theorem \ref{thm: tof intro} provides a new proof of the termination of pseudo-effective flips for lc pairs in dimension $\leq 3$ that does not rely on non-vanishing theorems. 

As an immediate consequence of Theorem \ref{thm: tof intro}, we obtain the following result:

\begin{thm}[Existence of log minimal models]\label{thm: emm intro}
    Let $(X, B,\Mm)/U$ be an NQC lc generalized pair of dimension $\leq 3$ over an algebraically closed field of characteristic $p>5$ such that $K_X+B+\Mm_X$ is pseudo-effective$/U$. Then:
    \begin{enumerate}
        \item $(X,B,\Mm)/U$ has a log minimal model (see Definition \ref{defn: model}).
        \item If $X$ is $\Qq$-factorial klt, then we may run a $(K_X+B+\Mm_X)$-MMP$/U$ and obtain a minimal model of $(X,B,\Mm)/U$.
    \end{enumerate}
\end{thm}
\begin{rem}
    In the case where $X$ is $\Qq$-factorial klt, we run the MMP by perturbing the original pair by an ample divisor; see Corollary \ref{cor: run MMP when X klt}.
    Following the ideas in \cite{HNT20}, we believe that one could try to 
    show the MMP can be run for essentially lc generalized pairs. We leave this question for future work.
\end{rem}

A key ingredient of the proof of Theorem \ref{thm: tof intro} is the ACC for lc thresholds for generalized pairs in dimension $\leq 3$.

\begin{thm}[ACC for lc thresholds]\label{thm: acc lct}
    Let $\Ii\subset [0,+\infty)$ be a DCC set. Then there exists an ACC set $\Ii'$, depending only on $\Ii$, with the following property: 
    Let $(X,B,\Mm)/U$ be an lc generalized pair of dimension $\leq 3$ over an algebraically closed field of characteristic $p>5$, $D\geq 0$ an $\Rr$-divisor on $X$, and $\Nn$ a nef$/U$ $\bb$-divisor on $X$ satisfying the following.
    Assume that
    \begin{enumerate}
        \item the coefficients of $B$ and $D$ belong to $\Ii$,
        \item $\Mm=\sum m_i\Mm_i$ and $\Nn=\sum n_i\Nn_i$, where each $\Mm_i,\Nn_i$ is a nef$/U$ $\bb$-Cartier $\bb$-divisor and $m_i,n_i\in\Ii$, and
        \item $D+\Nn_X$ is $\Rr$-Cartier.
    \end{enumerate}
    Then $$\lct(X,B,\Mm;D,\Nn):=\sup\{t\geq 0\mid (X,B+tD,\Mm+t\Nn)\text{ is lc}\}$$
    belongs to $\Ii'$.
\end{thm}

Another key ingredient of our proof is the existence of Shokurov-type polytopes:

\begin{thm}\label{thm: shokurov polytope intro}
    Let $(X,B=\sum_{j=1}^m b_jB_j,\Mm=\sum_{k=1}^n m_k\Mm_k)/U$ be an lc generalized pair of dimension $\leq 3$ over an algebraically closed field of characteristic $p>5$, where $B_j$ are the irreducible components of $B$, $m_k\geq 0$ for each $k$, and $\Mm_k$ is a nef$/U$ $\bb$-Cartier $\bb$-divisor for each $k$. Let $\bm{v}_0:=(b_1,\dots,b_m,m_1,\dots,m_n)\in\mathbb R^{m+n}$. 

    Assume that $K_X+B+\Mm_X$ is nef$/U$. Then there exists an open subset $V_0\ni\bm{v}_0$ of the rational envelope of $\bm{v}_0$ in $\mathbb R^{m+n}$, such that for any $(b_1',\dots,b_m',m_1',\dots,m_n')\in V_0$, $K_X+\sum_{j=1}^mb_j'B_j+\sum_{k=1}^nm_k'\Mm_{k,X}$ is nef$/U$ and $(X,\sum_{j=1}^mb_j'B_j,\sum_{k=1}^nm_k'\Mm_{k})$ is lc.
\end{thm}

\noindent\textbf{Acknowledgement}. This project is completed as a part of the research program ``Algebra and Number Theory Summer School Seminars" during the period of July 15th -- August 24th, 2024, at Peking University. The authors would like to thank the organizers and Peking University for their hospitality. They are especially grateful to their summer-school advisor, Jihao Liu, for his guidance, useful discussions, and comments. They also thank Liang Xiao and Shou-Wu Zhang for useful comments.

\section{Preliminaries}

We follow the standard notations and definitions as in \cite{BCHM10, KM98}. We will freely use the fact that log resolution exists for threefolds over an algebraically closed field \cite{CP18, CP19, Cut09}.

\subsection{Sets}

\begin{defn}
    Let $\Ii\subset [0,+\infty)$ be a set. We define $\Ii_+:=(\{\sum \gamma_i\mid \gamma_i\in\Ii\}\cup\{0\})\cap [0,1]$ and define $D(\Ii):=\{\frac{m-1+\gamma}{m}\mid m\in\mathbb N^+,\gamma\in\Ii_+\}$. We say that $\Ii$ satisfies the descending chain condition (DCC) if any decreasing sequence in $\Ii$ stabilizes. We say that $\Ii$ satisfies the ascending chain condition (ACC) if any non-decreasing sequence in $\Ii$ stabilizes.
\end{defn}

\subsection{\texorpdfstring{$\bb$}{}-divisors and generalized pairs}

We refer the reader to \cite[Definition 2.4]{HL23} for the definition and notation of $\bb$-divisors.

\begin{defn}
    Let $X$ be a normal quasi-projective variety of dimension $\leq 3$, $\tilde S$ be a prime divisor on $X$ with normalization $S$, and $\Dd$ be a $\bb$-divisor on $X$ which descends to a birational model of $X$. We define the \emph{restricted $\bb$-divisor} $\Dd|_S$ on $S$ in the following way. Let $h: X'\rightarrow X$ be a log resolution of $(X,\tilde S)$ such that $\Dd$ descends to $X'$ and let $S':=h^{-1}_*\tilde S$. We let $\Dd|_S:=\overline{\Dd_{X'}|_{S'}}$.
\end{defn}

\begin{defn}
    A \emph{generalized pair} $(X, B,\Mm)/U$ consists of a projective morphism $X\rightarrow U$ from a normal quasi-projective variety to a quasi-projective variety, an $\Rr$-divisor $B\geq 0$ on $X$, and a nef$/U$ $\bb$-divisor $\Mm$ on $X$, such that $K_X+B+\Mm_X$ is $\Rr$-Cartier. We say that $(X,B,\Mm)/U$ is NQC if $\Mm$ is NQC$/U$, i.e. $\Mm=\sum m_i\Mm_i$ where each $m_i\geq 0$ and each $\Mm_i$ is a nef$/U$ $\bb$-Cartier $\bb$-divisor. If $U$ is not important, then we may drop $U$. If $U=\{pt\}$, then we may drop $U$ and say that $(X, B,\Mm)$ is projective.

    For any prime divisor $E$ over $X$ and projective birational morphism $f: Y\rightarrow X$ such that $E$ is on $Y$ with $K_Y+B_Y+\Mm_Y:=f^*(K_X+B+\Mm_X)$, we define
    $$a(E,X,B,\Mm):=1-\mult_EB_Y$$
    as the \emph{log discrepancy} of $E$ with respect to $(X,B,\Mm)$. We say that $(X,B,\Mm)$ is \emph{klt} (resp. \emph{lc}) if $a(E,X,B,\Mm)>0$ (resp. $\geq 0$) for any prime divisor $E$ over $X$. We say that $E$ is an \emph{nklt place} of $(X,B,\Mm)$ if $a(E,X,B,\Mm)\leq 0$.

    An \emph{nklt center} of $(X, B,\Mm)$ is the center of an nklt place of $(X, B,\Mm)$ on $X$. The \emph{nklt locus} of $(X, B,\Mm)$, denoted by $\Nklt(X, B,\Mm)$, is the union of all nklt centers associated with a reduced subscheme structure. If $(X, B,\Mm)$ is lc, nklt places (resp. nklt centers) are also called \emph{lc places} (resp. \emph{lc centers}). We say that $(X,B,\Mm)$ is \emph{dlt} if $(X,B,\Mm)$ is lc, and for any lc center $W$ of $(X,B,\Mm)$ with generic point $\eta$, $\Mm$ descends to a neighborhood of $W$ and $(X,B)$ is log smooth near $\eta$.
\end{defn}

\begin{lem}\label{lem: hl22 3.4 char p}
    Let $(X,B,\Mm)/U$ be an lc generalized pair and $A$ an ample$/U$ $\Rr$-divisor on $X$. Assume that there exists a klt generalized pair $(X, B_0,\Mm_0)/U$. Then there exists a klt pair $(X,\Delta)$ such that
    $$\Delta\sim_{\mathbb R,U}B+\Mm_X+A.$$
\end{lem}
\begin{proof}
    Let $0<\epsilon\ll 1$ be a real number such that $A':=A+\epsilon(B+\Mm_X)-\epsilon(B_0+\Mm_{0,X})$ is ample$/U$. Possibly replacing $A$ with $A'$ and replacing $(X,B,\Mm)$ with $(X,(1-\epsilon)B+\epsilon B_0,(1-\epsilon)\Mm+\epsilon\Mm_0)$, we may assume that $(X,B,\Mm)$ is klt.

    By \cite[Theorem 2.13]{BMPSTWW23} there exists a log resolution $g: W\rightarrow X$ of $(X,\Supp B)$ such that $\Mm$ descends to $W$ and there exists a $g$-anti-ample $g$-exceptional divisor $E\geq 0$ on $W$. We write
    $$K_W+B_W+\Mm_W:=g^*(K_X+B+\Mm_X).$$
    Let $0<e\ll 1$ be a real number such that $\lfloor B_W+eE\rfloor=0$. Then $g^*A-eE$ is ample$/U$, so $\Mm_W+g^*A-eE$ is ample$/U$. Thus there exists $H_W\in |\Mm_W+g^*A-eE|_{\mathbb R}$ such that $(W,\Delta_W:=B_W+eE+H_W)$ is sub-klt. $\Delta:=g_*\Delta_W$ satisfies our requirements.
\end{proof}

As a corollary, we can run MMP for generalized pairs.

\begin{cor}\label{cor: run MMP when X klt}
    Let $(X,B,\Mm)/U$ be an lc generalized pair and $A$ an ample$/U$ $\Rr$-divisor on $X$. Assume that there exists a klt generalized pair $(X, B_0,\Mm_0)/U$.
    Then we can run a $(K_X+B+\Mm_X)$-MMP with scaling of $A$.
\end{cor}
\begin{proof}
    Since we can run MMP for the usual klt pair $(X,\Delta)$ (see \cite[Theorem 1.1]{HNT20}), the corollary follows from Lemma \ref{lem: hl22 3.4 char p}.
\end{proof}
\begin{rem}
    However, the termination of MMP for the usual pair does not imply the termination for the generalized pair since we need to perturb our pair by an ample divisor.
\end{rem}

\begin{lem}\label{lem: qfactorial dlt preserved under mmp}
    Let $(X,B,\Mm)/U$ be a $\Qq$-factorial lc generalized pair such that $X$ is klt, and let $(X,B,\Mm)\dashrightarrow (X',B',\Mm)$ be a sequence of steps of a $(K_X+B+\Mm_X)$-MMP$/U$. Then $X'$ is $\Qq$-factorial klt. Moreover, if $(X,B,\Mm)/U$ is dlt, then $(X',B',\Mm)/U$ is $\Qq$-factorial dlt.
\end{lem}
\begin{proof}
    We may assume that $\phi: X\dashrightarrow X'$ is a single step of a $(K_X+B+\Mm_X)$-MMP$/U$. 

    First, we show that $X'$ is $\Qq$-factorial klt. Then $X\dashrightarrow X'$ is a step of a $(K_X+B+\Mm_X+A)$-MMP$/U$ for some ample$/U$ $\Rr$-divisor $A$. By Lemma \ref{lem: hl22 3.4 char p}, there exists a klt pair $(X,\Delta)$ such that $\Delta\sim_{\mathbb R,U}B+\Mm_X+A$, so $X\dashrightarrow X'$ is a step of a $(K_X+\Delta)$-MMP$/U$. If $X\dashrightarrow X'$ is a divisorial contraction then the $\Qq$-factoriality follows from \cite[Lemma 7.2]{Bir16}. If $X\dashrightarrow X'$ is a flip with flipping contraction $X\rightarrow T$ and flipped contraction $X'\rightarrow T$, then $X'$ is the log minimal model of $(X,\Delta)$ over $T$ (\cite[Page 196, Line -2]{Bir16}), so $X'$ is $\Qq$-factorial. Moreover, let $\Delta'$ be the image of $\Delta$ on $X'$, then $(X',\Delta')$ is klt. Thus $X'$ is $\Qq$-factorial klt.

    Next, we show that $(X', B',\Mm)$ is dlt. Let $W'$ be an lc place of $(X',B',\Mm)$ with generic point $\eta'$ and let $E$ be a nklt place of $(X',B',\Mm)$ such that $\Center_{X'}E=W'$. Then
    $$0=a(E,X',B',\Mm)\geq a(E,X,B,\Mm)\geq 0,$$
    hence $a(E,X,B,\Mm)=0$ and $a(E,X',B',\Mm)=a(E,X,B,\Mm)$. Thus $\phi^{-1}$ is an isomorphism near $\eta'$. Let $W:=\Center_XE$ and let $\eta$ be the generic point of $W$. Then $W$ is an lc center of $(X,B,\Mm)$, so $(X,B)$ is log smooth near $\eta$ and  $\Mm$ descends over a neighborhood of $\eta$. Since $\phi^{-1}$ is an isomorphism near $\eta'$, $(X',B')$ is log smooth near $\eta'$ and $\Mm$ descends over a neighborhood of $\eta'$. Therefore, $(X',B',\Mm)$ is dlt.
\end{proof}

\begin{lem}\label{lem: dlt implies normal}
    Let $(X,B,\Mm)/U$ be a $\Qq$-factorial dlt generalized pair and $S$ an irreducible component of $\lfloor B\rfloor$. Then $S$ is normal.
\end{lem}
\begin{proof}
    By definition, $(X,B)$ is $\Qq$-factorial dlt, so the lemma follows from \cite[Lemma 5.2]{Bir16}
\end{proof}


\begin{defn}
    Let $(X, B,\Mm)/U$ be a generalized pair. Let $D\geq 0$ be an $\Rr$-divisor and $\Nn$ a nef$/U$ $\bb$-divisor such that $D+\Nn_X$ is $\Rr$-Cartier. We define
    $$\lct(X,B,\Mm;D,\Nn):=\sup\{t\geq 0\mid (X,B+tD,\Mm+t\Nn)\text{ is lc}\}$$
    to be the \emph{lc threshold} of $(D,\Nn)$ with respect to $(X,B,\Mm)$.
\end{defn}

\subsection{Models of generalized pairs}\label{sec:g-model}

\begin{defn}\label{defn: model}
    Let $(X,B,\Mm)/U$ and $(X',B',\Mm)/U$ be two lc generalized pairs associated with a birational map $\phi: X\dashrightarrow X'$ over $U$. We say that $(X',B',\Mm)/U$ is a \emph{log birational model} of $(X,B,\Mm)$ if $B'=\phi_*B+\Exc(\phi^{-1})$. 

    Assume that  $(X',B',\Mm)/U$ is a log birational model of $(X,B,\Mm)$ and $K_{X'}+B'+\Mm_{X'}$ is nef$/U$. We say that  $(X',B',\Mm)/U$ is a \emph{weak lc model} (resp. \emph{minimal model}) of  $(X,B,\Mm)/U$ in the sense of Birkar-Shokurov, or a \emph{bs-weak lc model} (resp. \emph{bs-minimal model}) of  $(X,B,\Mm)/U$, if $a(E,X,B,\Mm)\leq a(E,X',B',\Mm)$ (resp. $a(E,X,B,\Mm)<a(E,X',B',\Mm)$) for any prime divisor $E$ on $X$ that is exceptional$/X'$. In addition, if $X\dashrightarrow X'$ does not extract any divisor, then we say that $(X',B',\Mm)/U$ is a \emph{weak lc model} (resp. \emph{minimal model}) of $(X,B,\Mm)/U$.
    
    We say that  $(X',B',\Mm)/U$ is a \emph{log minimal model} of  $(X,B,\Mm)/U$ if $(X',B',\Mm)/U$ is a bs-minimal model of $(X,B,\Mm)/U$ and $(X',B',\Mm)/U$ is $\Qq$-factorial dlt.
\end{defn}

\section{Technical preparations}

\subsection{Adjunction formula for generalized pairs}

\begin{lem}\label{lem: adjunction}
    Let $(X,B,\Mm)/U$ be an lc generalized pair of dimension $\leq 3$, $\tilde S$ an irreducible component of $\lfloor B\rfloor$, and $S$ the normalization of $\tilde S$.  Assume that $B=\tilde S+\sum_{j=1}^m b_jB_j$ and $\Mm=\sum_{k=1}^n r_k\Mm_k$, where $B_j$ are the irreducible components of $B$, and each $\Mm_k$ is a nef$/U$ $\bb$-Cartier $\bb$-divisor. Let $\Mm_k^S:=\Mm_k|_S$ for each $k$ and let $\Mm^S:=\Mm|_S$. Then there exist prime divisors $T_1,\dots,T_l,C_1,\dots,C_q$ on $S$, positive integers $w_1,\dots,w_q$, and non-negative integers $\{d_{i,j}\}_{1\leq i\leq q, 1\leq j\leq m}$ and $\{e_{i,k}\}_{1\leq i\leq q,1\leq k\leq n}$, such that for any real numbers $b_1',\dots,b_m',r_1',\dots,r_n'$, we have the following. Let $B':=\tilde S+\sum_{j=1}^mb_j'B_j$ and $\Mm':=\sum_{k=1}^nr_k'\Mm_k$. Assume that $K_{X}+B'+\Mm'_{X}$ is $\Rr$-Cartier. Then:
    \begin{enumerate}
        \item 
        $$K_{S}+B'_{S}+\Mm'^S_{S}:=(K_X+B'+\Mm'_X)|_S,$$
        where $\Mm'^S:=\Mm'|_S$, and
        $$B'_{S}=\sum_{i=1}^lT_i+\sum_{i=1}^q\frac{w_i-1+\sum_{j=1}^md_{i,j}b_j'+\sum_{k=1}^ne_{i,k}r_k'}{w_i}C_i.$$
        In particular, if $b_j'$ and $r_k'$ belong to a DCC set, then the coefficients of $B_S'$ belong to a DCC set.
        \item If $(X,B',\Mm')$ is lc, then $(S,B'_S,\Mm'^S)$ is lc.
    \end{enumerate}
\end{lem}
\begin{proof}
    By \cite[Remark 2.5]{FW23}, we have (2), so we only need to prove (1). We only need to determine the coefficient of $V$ in $B'_S$ for any irreducible component $V$ of $B'_S$. Fix such a component $V$ and let $W$ be its image in $X$. Determining $\mult_VB_S'$ is a local problem near the generic point of $W$, so by replacing $X$ with $\Spec \mathcal{O}_{X, W}$, we may assume that $X$ is a normal excellent scheme of dimension $2$, $\dim S=1$, and $V$ is a closed point.

    Since $(X,B,\Mm)$ is lc, $(X,\tilde S)$ is numerically lc (cf. \cite[Notation 4.1]{KM98}). If $(X,\tilde S)$ is not numerically plt near $V$, then $B=\tilde S$ near $V$ and $\Mm$ descends to $X$ near a neighborhood of $V$, and $\mult_VB_S=1$. Thus $\mult_VB'_S=1$ and we are done. Therefore, we may assume that $(X,\tilde S)$ is numerically plt near $V$, hence $X$ is numerically klt near the generic point of $V$. Thus $X$ is $\Qq$-factorial klt near $V$. Possibly shrinking $X$ to a neighborhood of $V$, we may assume that $X$ is $\Qq$-factorial.

    By \cite[3.35]{Kol13}, there exists a positive integer $w=w_V$, such that for any Weil divisor $D$ on $X$, $wD$ is Cartier near $V$. We have $B_j|_S:=\frac{d_j}{w}V$ near a neighborhood of $V$ for some non-negative integers $d_j$. 

    Let $f: X'\to X$ be a log resolution of $(X,B)$ such that $\Mm_k$ descends to $X'$ for each $k$. By the negativity lemma,
    $$f^*\Mm_{k,X}=\Mm_{k,X'}+E_k$$
    for some $\Qq$-divisors $E_k\geq 0$. Since $\Mm_{k,X'}$ is Cartier and $w\Mm_{k,X}$ is Cartier near $V$, $wE_k$ is Cartier over a neighborhood of $V$. Let $S':=f^{-1}_*\tilde S$ and let $f_S: S'\rightarrow S$ be the induced birational morphism, then
    $$\Mm_{k,X}|_S=(f_*(\Mm_{k,X'}+E_k))|_S=(f_S)_*(\Mm^S_{k,S'}+E_k|_{S'})=\Mm^S_{k,S}+(f_S)_*E_k|_{S'},$$
    Let $e_k:=w\mult_V((f_S)_*E_k|_{S'})$, then $e_k\in\mathbb N$. 
    By \cite[Proposition 4.2]{Bir16}, 
    $$K_S+\frac{w-1}{w}V=(K_X+S)|_S$$
    near $V$. Thus
    $$\mult_VB'_S=\frac{w-1}{w}+\sum b'_j\frac{d_j}{w}+\sum r_k'\frac{e_k}{w}.$$
    (1) follows.
\end{proof}

\begin{lem}\label{lem: dlt under adjunction}
    Let $(X, B,\Mm)/U$ be a $\Qq$-factorial dlt generalized pair of dimension $\leq 3$ and let $S$ be a component of $\lfloor B\rfloor$. Let $(S,B_S,\Mm^S)/U$ be the generalized pair induced by adjunction
    $$K_S+B_S+\Mm^S_S:=(K_X+B+\Mm_X)|_S.$$
    Then:
    \begin{enumerate}
        \item $(S,B_S,\Mm^S)$ is dlt.
        \item For any lc center $V\subset S$ of $(X,B,\Mm)$, $V$ is an lc center of $(S,B_S,\Mm^S)$.
    \end{enumerate}
\end{lem}
\begin{proof}
    By Lemma \ref{lem: dlt implies normal}, $S$ is normal. Let $h: X'\rightarrow X$ be a log resolution of $(X,B)$ such that $\Mm$ descends to $X'$, $K_{X'}+B'+\Mm_{X'}:=h^*(K_X+B+\Mm_X)$, and $S':=h^{-1}_*S$. Let
    $$K_{S'}+B_{S'}+\Mm^S_{S'}:=(K_{X'}+B'+\Mm_{X'})|_{S'}$$
    where $B_{S'}:=(B'-S')|_{S'}$, and let $h_S: S'\rightarrow S$ be the induced birational morphism. Then
    $$K_{S'}+B_{S'}+\Mm^S_{S'}=h_S^*(K_S+B_S+\Mm^S_S).$$
    For any lc place $E_S$ of $(S,B_S,\Mm^S)$, $\Center_{S'}E_S$ is a stratum of $B_{S'}^{=1}$, hence a stratum of $B'^{=1}$. Let $E$ be an lc place of $(X,B,\Mm)$ such that $\Center_{X'}E=\Center_{S'}E_S$. Since $(X, B,\Mm)$ is dlt, $(X, B)$ is log smooth near the generic point of $\Center_XE$ and $\Mm$ descends over a neighborhood of the generic point of $\Center_XE$. Thus $(S, B_S)$ is log smooth near the generic point of $\Center_XE=\Center_{S}E_S$ and $\Mm^S$ descends over a neighborhood of the generic point of $\Center_SE_S$. Therefore, $(S,B_S,\Mm^S)$ is dlt.

    For any lc center $V$ of $(X,B,\Mm)$, possibly shrinking $X$ to a neighborhood of $V$, we have that $(X,B)$ is log smooth, $V$ is a stratum of $\lfloor B\rfloor$, and $\Mm$ descends to $X$. Thus $(S,B_S)$ is log smooth, $V$ is a stratum of $\lfloor B_S\rfloor$ $V$, and $\Mm^S$ descends to $S$. In particular, $V$ is an lc center of $(S,B_S,\Mm^S)$.
\end{proof}

\subsection{Special termination}

\begin{prop}\label{prop: special termination}
    Let $(X,B,\Mm)/U$ be a $\Qq$-factorial NQC dlt generalized pair of dimension $3$. Then for any sequence of $(K_X+B+\Mm_X)$-flips$/U$, after finitely many flips, the flipping locus does not intersect the strict transform of $\Supp\lfloor B\rfloor$.
\end{prop}
\begin{proof}
    The proof almost follows from the same lines of the proof of \cite[Proposition 2.16]{FW23}, which proves the case when $B$ is a $\Qq$-divisor and $\Mm$ is a $\Qq$-$\bb$-divisor. For the reader's convenience, we provide a full proof here.
    
    \noindent\textbf{Step 1}. We step up the MMP and deal with dimension $0$ lc centers in this step. Let
    $$(X,B,\Mm)=:(X_0,B_0,\Mm)\dashrightarrow (X_1,B_1,\Mm)\dashrightarrow (X_2,B_2,\Mm)\dashrightarrow\dots\dashrightarrow (X_i,B_i,\Mm)\dashrightarrow\dots$$
    be an infinite sequence of $(K_X+B+\Mm_X)$-flips$/U$. By Lemma \ref{lem: qfactorial dlt preserved under mmp}, $(X_i,B_i,\Mm)$ is $\Qq$-factorial dlt for any $i$. Possibly truncating the MMP, we may assume that no lc center is contracted by this MMP. In particular, the MMP is an isomorphism near the generic point of any lc center, and in particular, is an isomorphism near any dimension $0$ lc center.

    \smallskip

    \noindent\textbf{Step 2}. We deal with dimension $1$ lc centers in this step. Let $C$ be an lc center of $(X, B,\Mm)$ such that $\dim C=1$. Since $(X, B,\Mm)$ is dlt, there are two irreducible components $S, T$ of $\lfloor B\rfloor$ such that $C$ is a component of $S\cap T$. For each $i$, let $S_i,T_i$ be the strict transforms of $S,T$ on $X_i$, and let $C_i$ be the image of $C$ on $X_i$. By Lemma \ref{lem: dlt implies normal}, $S_i, T_i$ are normal. Let $(S_i,B_{S_i},\Mm^S)/U$ be the generalized pair induced by adjunction
    $$K_{S_i}+B_{S_i}+\Mm^{S_i}_{S_i}:= (K_{X_i}+B_{i}+\Mm_{X_i})|_{S_i}.$$ 
    Since the MMP is an isomorphism near the generic point of $S_i$ for each $i$, the $\bb$-divisor $\Mm^{S_i}$ does not depend on $i$, and we may denote it by $\Mm^S$.
    
    By Lemma \ref{lem: dlt under adjunction}, $(S_i,B_{S_i},\Mm^S)/U$ is dlt and $C_i$ is an lc center of $(S_i,B_{S_i},\Mm^S)$. Since $S_i$ is a surface, $(S_i,B_{S_i})$ is numerically dlt, hence $S_i$ is numerically klt. Thus $S_i$ is $\Qq$-factorial, so $(S_i,B_{S_i},\Mm^S)/U$ is $\Qq$-factorial dlt.  By Lemma \ref{lem: dlt implies normal}, $C_i$ is normal. Let $(C_i,B_{C_i},\Mm^C)/U$ be the generalized pair induced by adjunction
    $$K_{C_i}+B_{C_i}+\Mm^{C_i}_{C_i}:= (K_{S_i}+B_{S_i}+\Mm^S_{S_i})|_{C_i}.$$
    Since the MMP is an isomorphism near the generic point of $C_i$ for each $i$, the $\bb$-divisor $\Mm^{C_i}$ does not depend on $i$, and we may denote it by $\Mm^C$. Since $C_i$ are curves, $C_i\cong C$ for any $i$, hence $\Mm^C_{C_i}=\Mm^C_{C_{i+1}}$ for any $i$.
    
    Since $a(D,X_i,B_i,\Mm)\leq a(D,X_{i+1},B_{i+1},\Mm)$ for each $i$ and any prime divisor $D$ over $X$, by adjunction,
    $$B_{C_i}\geq B_{C_{i+1}}$$
    for each $i$. Here we identify $B_{C_i}$ and $B_{C_{i+1}}$ as $\Rr$-divisors on $C$.  By Lemma \ref{lem: adjunction}, there exists a DCC set $\Ii$ depending only on $B_0$ and $\Mm$ such that the coefficients of $B_{C_i}$ and $B_{C_{i+1}}$ belong to $\Ii$. Therefore, by truncating the MMP, we have $B_{C_i}=B_{C_{i+1}}$ for any $i$. Therefore,
    $$K_{C_i}+B_{C_i}+\Mm^C_{C_i}=K_{C_{i+1}}+B_{C_{i+1}}+\Mm^{C}_{C_{i+1}}$$
    for any $i$. Thus, the MMP is an isomorphism near $C$.

    \smallskip

    \noindent\textbf{Step 3}. We deal with dimension $2$ lc centers in this step and conclude the proof. 

    Let $S$ be an irreducible component of $\lfloor B\rfloor$. For each $i$, let $X_i\rightarrow T_i\leftarrow X_{i+1}$ be each step of the MMP and let $E_i,F_i$ be the exceptional locus of $X_i\rightarrow T_i$, $X_{i+1}\rightarrow T_i$.


    Suppose that $S_i$ contains a component of $F_{i-1}$ and let $D_i$ be such a component. Since $B_{S_i}\geq 0$, $a(D_i,S_0,B_{S_0},\Mm^S)<a(D_i,S_i,B_{S_i},\Mm^S)\leq 1$. If there exists an lc center of $(S_0,B_{S_0},\Mm^S)$ which contains $\Center_{S_0}D_i$, then there exists a component $T_0$ of $\lfloor B_0\rfloor$ such that $D_i$ is contained in $S_0\cap T_0$. This contradicts \textbf{Step 2}. Thus $\Center_{S_0}D_i$ is not contained in $\Nklt(S_0,B_{S_0},\Mm^S)$. Since $S_0$ is a surface, there are finitely many prime divisors over $S_0$ whose centers on $S_0$ are not contained in $\Nklt(S_0, B_{S_0},\Mm^S)$ and with log discrepancies $\leq 1$ with respect to $\Nklt(S_0, B_{S_0},\Mm^S)$. Thus, by identifying $D_i$ with its image on $X_j$ for any $j$, there are only finitely many different choices of $D_i$. Since the coefficient of $D_i$ in $B_{S_j}$ belongs to a DCC set, possibly truncating the MMP, we may assume for any $D_i$ and any $j$, the coefficient of $D_i$ in $B_{S_j}$ is a constant. Therefore, possibly truncating the MMP, we may assume that the $S_i$ does not contain any component of $F_{i-1}$.

    Finally, suppose that $S_i$ intersects $F_{i-1}$ in finitely many points. Then $S_i$ is anti-ample$/T_{i-1}$ and then $S_{i-1}$ is ample$/T_{i-1}$. Therefore, at least one component of $E_{i-1}$ is contained in $S_{i-1}$, hence the induced morphism $S_{i-1} \to S_i$ contracts at least a curve, so $\rho(S_i)<\rho(S_{i-1})$. This can happen only finitely many times. Thus, possibly truncating the MMP, $S_i$ does not intersect $F_i$ in finitely many points. 

    In summary, $S_i$ does not intersect $F_i$, so the MMP is an isomorphism near $S_i$, and we are done.
\end{proof}

\subsection{Existence of dlt modification}

\begin{lem}\label{lem:extract divisor}
    Let $(X, B,\Mm)$ be an NQC lc generalized pair of dimension $\leq 3$. Let $S_1,\dots,S_r$ be prime divisors that are exceptional$/X$ such that $a(S_i,X,B,\Mm)\leq 1$ for each $i$. Then there exists a projective birational morphism $f: Y\rightarrow X$ satisfying the following:
    \begin{enumerate}
        \item $S_1,\dots,S_r$ are divisors on $Y$, 
        \item for any $f$-exceptional prime divisor $E$, either $E=S_i$ for some $i$ or $a(E,X,B,\Mm)=0$.
        \item $(Y,B_Y,\Mm)$ is $\Qq$-factorial dlt, where $K_Y+B_Y+\Mm_Y:=f^*(K_X+B+\Mm_X)$.
    \end{enumerate}
\end{lem}
\begin{proof}
    It essentially follows from the same lines of the proof of \cite[Proposition 2.9]{BF25}, which proves the case when $B$ is a $\Qq$-divisor and $\Mm$ is a $\Qq$-$\bb$-divisor. For the reader's convenience, we provide a full proof here.

    Let $g: W\rightarrow X$ be a log resolution such that $\Mm$ descends to $W$, and let $E_1,\dots,E_n$ be all $g$-exceptional prime divisors, where $E_i=S_i$ for any $1\leq i\leq r$. Let $a_i:=a(E_i,X,B,\Mm)$ for any $i$ and let
    $$K_W+B_W+\Mm_W:=g^*(K_X+B+\Mm_X)+\sum_{j>r}a_jE_j.$$
    Then $(W,B_W,\Mm)$ is dlt and $K_W+B_W+\Mm_W\sim_{\mathbb R,X}\sum_{j>r} a_jE_j$.

    We run a $(K_W+B_W+\Mm_W)$-MMP$/X$. By special termination (Proposition \ref{prop: special termination}), this MMP terminates near the image of $\lfloor B_W\rfloor=\cup_{i=1}^nE_i$. Since the $(K_W+B_W+\Mm_W)$-MMP$/X$ is also a $(\sum_{i>r} a_iE_i)$-MMP$/X$, this MMP terminates with a minimal model $(Y,B_Y,\Mm)/X$ of $(W,B_W,\Mm)/X$, and any divisor contracted by this MMP is $E_i$ for some $i>r$ such that $a_i>0$. On the other hand, let $E_{i,Y}$ be the image of $E_i$ on $Y$ for each $i$, then $\sum_{i>r}a_iE_{i,Y}$ is nef$/X$, so by the negativity lemma, $\sum_{i>r}a_iE_{i,Y}=0$, hence
    $$K_Y+B_Y+\Mm_Y=f^*(K_X+B+\Mm_X),$$
    and for any $i>r$ such that $a_i>0$, $E_i$ is contracted by this MMP. Therefore, the divisors extracted by the induced birational morphism $f: Y\rightarrow X$ are exactly $E_1,\dots, E_r$ and the $E_i$ such that $i>r$ and $a_i=0$. Since $(W,B_W,\Mm)$ is $\Qq$-factorial dlt, so is $(Y,B_Y,\Mm)$. The lemma follows.
\end{proof}

\begin{cor}\label{cor: extistence of terminalization}
    Let $(X, B,\Mm)$ be a klt generalized pair. 
    Then there exists a birational morphism $f:Y \to X$ such that $(Y,B_Y,\Mm)$ is $\Qq$-factorial terminal with $K_Y+B_Y+\Mm_Y = f^*(K_X+B+\Mm_X)$.
    Such a generalized pair $(Y, B_Y,\Mm)$ together with the birational morphism $f: Y \to X$ is called a \emph{terminalization} of $(X, B,\Mm)$.
\end{cor}
\begin{proof}
    Since $(X, B,\Mm)$ is klt, there are only finitely many prime divisors $S_1,\ldots,S_r$ over $X$ with discrepancy $a(S_i,X,B,\Mm)\leq 1$.
    Let $(Y, B_Y,\Mm)$ be the generalized pair constructed in Lemma \ref{lem:extract divisor}.
    Then $(Y,B_Y,\Mm)$ is $\Qq$-factorial dlt such that $K_Y+B_Y+\Mm_Y = f^*(K_X+B+\Mm_X)$.
    For any other prime divisor $E$ which is exceptional over $Y$, it is also exceptional over $X$.
    By our construction, either $E = S_i$ for some $i$ or $a(E,X,B,\Mm) > 1$.
    Since $E$ is exceptional over $Y$, $E \neq S_i$ and hence $a(E,Y,B_Y,\Mm) = a(E,X,B,\Mm) > 1$.
    Hence $(Y,B_Y,\Mm)$ is terminal.
\end{proof}

\begin{deflem}\label{deflem: dlt model}
    Let $(X, B,\Mm)$ be an lc generalized pair. A \emph{$\Qq$-factorial dlt modification} of $(X,B,\Mm)$ is a birational morphism $f: Y\rightarrow X$ such that
    \begin{enumerate}
        \item For any $f$-exceptional prime divisor $E$, $a(E,X,B,\Mm)=0$.
        \item $(Y,B_Y,\Mm)$ is $\Qq$-factorial dlt, where $K_Y+B_Y+\Mm_Y:=f^*(K_X+B+\Mm_X)$.
    \end{enumerate}
    We say that $(Y,B_Y,\Mm)$ is a \emph{dlt model} of $(X,B,\Mm)$. By Lemma \ref{lem:extract divisor}, for NQC lc generalized pairs of dimension $\leq 3$, dlt models exist.
\end{deflem}

\begin{lem}\label{lem: special extraction}
    Let $(X,B,\Mm)$ be an NQC lc generalized pair of dimension $\leq 3$, $D\geq 0$ an $\Rr$-divisor on $X$, and $\Nn$ a nef$/X$ $\bb$-divisor on $X$ such that $D+\Nn_X$ is $\Rr$-Cartier. Assume that 
    \begin{itemize}
        \item[(i)] $(X,B+D,\Mm+\Nn)$ is lc,
        \item[(ii)]  $(X,B+(1+\epsilon)D,\Mm+(1+\epsilon)\Nn)$ is not lc for any $\epsilon>0$, and
        \item[(iii)] for any prime divisor $P$ on $X$ such that $\mult_P(B+D)=1$, $\mult_PD=0$.
    \end{itemize}
    Then for any real number $t\in (0,1)$, there are two projective birational morphisms $h: X'\rightarrow X$, $g: Y\rightarrow X'$ satisfying the following.
    \begin{enumerate}
        \item $h$ is a $\Qq$-factorial dlt modification of $(X,B+tD,\Mm+t\Nn)$.
        \item Any $h$-exceptional prime divisor $P$ is an lc place of $(X,B,\Mm)$. In particular, $\mult_P(D+\Nn_X)=0$ and $a(P,X,B+sD,\Mm+s\Nn)=0$ for any real number $s$.
        \item $g$ extracts a unique prime divisor $E$. In particular, $-E$ is ample$/X'$.
        \item $a(E,X,B+D,\Mm+\Nn)=0$ and $a(E,X,B,\Mm)>0$.
        \item Let $F$ be the reduced $(h\circ g)$-exceptional divisor and $B_Y,D_Y$ the strict transforms of $B,D$ on $Y$ respectively. Then $(Y,B_Y+tD_Y+F,\Mm+t\Nn)$ is $\Qq$-factorial dlt.
    \end{enumerate}
\end{lem}
\begin{proof}
    By condition (ii), there exists a prime divisor $P_0$ over $X$ such that $a(P_0, X, B+D,\Mm+\Nn)=0$ and $a(P_0, X, B+(1+\epsilon)D,\Mm+(1+\epsilon)\Nn)<0$ for any $\epsilon>0$. By condition (iii), $P_0$ is exceptional$/X$. By Lemma \ref{lem:extract divisor}, there exists a $\Qq$-factorial dlt modification $f: W\rightarrow X$ of $(X,B+D,\Mm+\Nn)$ such that $P_0$ is on $W$. Let $B_W,D_W$ be the strict transforms of $B$ and $D$ on $W$ respectively, $F_1,\dots,F_n$ be the prime $f$-exceptional divisors, and $F_W=\sum_{i=1}^nF_i$ the reduced $f$-exceptional divisor. Let $a_i:=a(F_i,X,B+tD,\Mm+t\Nn)$. Then
    $$K_W+B_W+D_W+F_W+\Mm_W+\Nn_W=f^*(K_X+B+D+\Mm_X+\Nn_X).$$
    Moreover, $(W,B_W+tD_W+F_W,\Mm+t\Nn)$ is dlt, and we have
    $$K_W+B_W+tD_W+F_W+\Mm_W+t\Nn_W\sim_{\mathbb R,X}\sum_{i=1}^na_iF_i.$$
    We run a $(K_W+B_W+tD_W+F_W+\Mm_W+t\Nn_W)$-MMP$/X$. By special termination (Proposition \ref{prop: special termination}), this MMP terminates near the images of $F_W$. Since this MMP is also a $(\sum_{i=1}^na_iF_i)$-MMP$/X$, this MMP terminates with a minimal model $(X',B'+tD'+F',\Mm+t\Nn)/X$ of $(W,B_W+tD_W+F_W,\Mm+t\Nn)/X$, where $B',D',F'$ are the images of $B_W,D_W,F_W$ on $X'$ respectively. Let $F_{i,W}$ be the image of $F_i$ on $W$ for each $i$, then $\sum a_iF_{i,W}$ is nef$/X$. By the negativity lemma, $F_i$ is contracted by the MMP if $a_i>0$. Moreover, since the MMP is a $(\sum_{i=1}^na_iF_i)$-MMP, the MMP only contracts divisors that are contained in $F_W$. Therefore, the divisors contracted by this MMP are exactly the $F_i$ such that $a_i>0$. In particular, $(X',B'+tD'+F',\Mm+t\Nn)$ is a dlt model of $(X,B+tD,\Mm+t\Nn)$, and $P_0$ is contracted by this MMP. We let $h: X'\rightarrow X$ be the induced morphism. Then (1) holds.

    Since the induced birational map $W\dashrightarrow X'$ contracts at least one divisor, $W\dashrightarrow X'$ is not the identity morphism. We let $g: Y\dashrightarrow X'$ be the last step of this MMP. Since $X'$ is $\Qq$-factorial and
    $$K_{X'}+B'+tD'+F'+\Mm_{X'}+t\Nn_{X'}\sim_{\mathbb R,X}0,$$
    $g$ is a divisorial contraction of a divisor $E$.

    We show that $g$ and $h$ satisfy the required properties. We have already shown (1), and (5) follows from our construction. (3) follows from our construction and the negativity lemma. For any $h$-exceptional prime divisor $P$, we have
    $$a(P,X,B+tD,\Mm+t\Nn)=0=a(P,X,B+D,\Mm+\Nn)$$
    as $P$ is also extracted by $f$. This implies (2). Finally, since $E$ is extracted by $f$, $a(E,X,B+D,\Mm+\Nn)=0$. Let $B_Y,D_Y,F_Y$ be the images of $B_W,D_W,F_W$ on $Y$ respectively. Since $g$ is a $(K_Y+B_Y+tD_Y+F_Y+\Mm_Y+t\Nn_Y)$-negative contraction, by (2), we have
    $$a(E,X,B+tD,\Mm+t\Nn)=a(E,X',B'+tD'+F',\Mm+t\Nn)>a(E,Y,B_Y+tD_Y+F_Y,\Mm+t\Nn)\geq 0.$$
    This implies (4), and we are done.
\end{proof}

\section{ACC for lc thresholds}

\begin{nota}
    Let $\Ii\subset\mathbb R$ be a set. Let $X\rightarrow U$ be a projective morphism from a normal quasi-projective variety to a quasi-projective variety. For any $\Rr$-divisor $D$ on $X$, we write $D\in\Ii$ by abuse of notation if all coefficients of $D$ belong to $\Ii$. For any NQC$/U$ $\bb$-divisor $\Mm$ on $X$, we write $\Mm\in\Nef(U,\Ii)$ if $\Mm=\sum m_i\Mm_i$, where each $\Mm_i$ is a nef$/U$ $\bb$-Cartier $\bb$-divisor and each $m_i \in \Gamma$. If, in addition, $\Mm_i\not\equiv_U\bm{0}$ for any $i$, then we may write $\Mm\in\Nef^0(U,\Ii)$. If $U=\{pt\}$ then we may drop $U$ and use $\Nef(\Ii)$ and $\Nef^0(\Ii)$ respectively.
\end{nota}

\begin{conj}\label{conj: ACC for lc thresholds}
    Let $d$ be a positive integer and $\Ii\subset [0,+\infty)$ a DCC set. Then
    $$\{\lct(X,B,\Mm;D,\Nn)\mid \dim X=d,B,D\in\Ii,\Mm,\Nn\in\Nef(X,\Ii)\}$$
    satisfies the ACC.
\end{conj}

\begin{conj}\label{conj: global ACC}
    Let $d$ be a positive integer and $\Ii\subset [0,+\infty)$ a DCC set. Then there
    exists a finite set $\Ii_0\subset\Ii$ depending only on $d$ and $\Ii$ satisfying the following. Assume that
    \begin{enumerate}
        \item $(X,B,\Mm)$ is a projective lc generalized pair of dimension $d$,
        \item $B\in\Ii$ and $\Mm\in\Nef^0(\Ii)$, and
        \item $K_X+B+\Mm_X\equiv 0$.
    \end{enumerate}
    Then $B\in\Ii_0$ and $\Mm\in\Nef^0(\Ii_0)$.
\end{conj}

\subsection{ACC for lc thresholds on surfaces}

\begin{lem}\label{lem:surface ACC}
Conjecture \ref{conj: ACC for lc thresholds} holds when $d=2$.
\end{lem}
\begin{proof}
We may assume that $1\in\Ii$. Suppose the statement is false, then there exists a sequence $(X_i, B_i,\Mm_i; D_i,\Nn_i)$ such that $\dim X_i=2$, $B_i, D_i\in\Ii$, $\Mm,\Nn\in\Nef(X,\Ii)$, and $t_i:=\lct(X_i, B_i,\Mm_i;D_i,\Nn_i)$ is strictly increasing. In particular, we may assume that $t_i>0$ for each $i$. 

    Suppose that for infinitely many $i$, there exists an lc place $E_i$ of $(X_i,\Delta_i:=B_i+t_iD_i,\Pp_i:=\Mm_i+t_i\Nn_i)$ on $X_i$ such that
    $$a(E_i,X_i,B_i+(t_i+\epsilon)D_i,\Mm_i+(t_i+\epsilon)\Nn_i)<0$$
    for any $\epsilon>0$. Then $E_i$ is a component of $B_i+t_iD_i$, $\mult_{E_i}D_i>0$, and $\mult_{E_i}(B_i+t_iD_i)=1$, which is not possible as $\Ii$ satisfies the DCC and $t_i$ is strictly increasing. Thus, possibly passing to a subsequence, we may assume that no such $E_i$ exists. 
    
    By Lemma \ref{lem: special extraction}, there are two projective birational morphisms $h_i: X_i'\rightarrow X_i$ and $g_i: Y_i\rightarrow X_i'$, such that
\begin{itemize}
    \item   $h_i$ is a $\Qq$-factorial dlt modification of $(X_i,\Delta_i,\Pp_i)$,
    \item any $h_i$-exceptional prime divisor $P_i$ is an lc place of $(X_i,B_i,\Mm_i)$,
    \item $g_i$ extracts a unique prime divisor $E_i$ such that $-E_i$ is ample$/X_i'$, 
    \item $a(E_i,X_i,B_i+t_iD_i,\Mm_i+t_i\Nn_i)=0$ and $a(E_i,X_i,B_i,\Mm_i)>0$, and
    \item $(Y_i,B_{Y_i},\Mm)$ is $\Qq$-factorial dlt, where $B_{Y_i}$ is the strict transform of $B_i$ on $Y_i$ plus the reduced $(h_i\circ g_i)$-exceptional divisor. By Lemma \ref{lem: dlt implies normal}, $E_i$ is normal.
\end{itemize}
Let $D_i',D_{Y_i}$ be the strict transforms of $D_i$ on $X_i'$, $Y_i$ respectively, and let $B_i'$ be the strict transform of $B_i$ on $X_i'$ plus the reduced $h_i$-exceptional divisor. By our construction,
$$t_i=\lct(X_i',B_i',\Mm_i;D_i',\Nn_i),$$
and $D_{Y_i}+\Nn_{i,Y_i}=g_i^*(D_i'+\Nn_{i,X_i'})+b_iE_i$ for some $b_i>0$. Therefore, $(D_{Y_i}+\Nn_{i,Y_i})|_{E_i}$ is ample. 

Let $x_i':=\Center_{X_i'}E_i$ and let
$$K_{E_i}+B_{E_i}(t)+\Mm_i(t)_{E_i}:=(K_{Y_i}+B_{Y_i}+tD_{Y_i}+\Mm_{i,Y_i}+t\Nn_{i,Y_i})|_{E_i}$$
for any real number $t$. Then $$K_{E_i}+B_{E_i}(t_i)+\Mm_i(t_i)_{E_i}\equiv 0$$
and
 $$K_{E_i}+B_{E_i}(0)+\Mm_i(0)_{E_i}$$
 is anti-ample. Since $B_{E_i}(0)+\Mm_i(0)_{E_i}$ is pseudo-effective, $E_i\cong\mathbb P^1$. Therefore, 
 $$2=\deg(B_{E_i}(t_i)+\Mm_i(t_i)_{E_i}).$$
By Lemma \ref{lem: adjunction}, the coefficients of $B_{E_i}(t_i)$ belong to a DCC set. Thus possibly passing to a subsequence, the coefficients of $B_{E_i}(t_i)$ belong to a finite set, and $\Nn_{i,Y_i}|_{E_i}\equiv 0$. Thus $D_{Y_i}|_{E_i}\not=0$. Since $t_i$ is strictly increasing, by Lemma \ref{lem: adjunction} again, the coefficients of $B_{E_i}(t_i)$ belong to a finite set, a contradiction.
\end{proof}

\subsection{Global ACC for surfaces}

\begin{lem}\label{lem:global ACC}
Conjecture \ref{conj: global ACC} holds when $d=2$.
\end{lem}
\begin{proof}
We may assume that $1\in\Ii$. Suppose that the conjecture does not hold, then there exists a sequence $(X_i, B_i:=\sum b_{i,j}B_{i,j},\Mm_i:=\sum m_{i,j}\Mm_{i,j})$ as in the assumptions, such that $b_{i,j},m_{i,j}\in\Ii$, either $b_{i,1}$ is strictly increasing or $m_{i,1}$ is strictly increasing, $B_{i,j}$ are distinct prime divisors for any fixed $i$, and $\Mm_{i,j}\not\equiv\bm{0}$ for any $i,j$.  By Definition-Lemma \ref{deflem: dlt model}, possibly replacing $(X_i,B_i,\Mm_i)$ with a dlt model, we may assume that $X_i$ is $\Qq$-factorial klt.

\medskip

\noindent\textbf{Step 1}. We prove the statement when there exists a contraction $X_i\rightarrow Z_i$ such that $Z_i$ is a curve, $B_{i,1}$ is horizontal$/Z_i$ if $b_{i,1}$ is strictly increasing, and $\Mm_{i,1,X_i}$ is ample$/Z_i$ if $b_{i,1}$ is not strictly increasing.

Let $F_i$ be the reduced variety associated to a general fiber of $X_i\rightarrow Z_i$. Similarly to the proof of \cite[Proposition 11.7]{Bir16}, we can make the following computation:
$$(K_{X_i}+F_i)\cdot F_i = 2p_a(F_i)-2, \quad F_i^2=0.$$
Hence, we have
$$0=(K_{X_i}+B_i+\Mm_{i,X_i})\cdot F_i=2p_a(F_i)-2+\sum b_{i,j}(B_{i,j}\cdot F_i)+\sum m_{i,j}(\Mm_{i,j,X_i}\cdot F_i),$$
and $B_{i,j}\cdot F_i,\Mm_{i,j,X_i}\cdot F_i\in\mathbb N$. 
So there exists a finite set $\Ii_0\subset\Ii$ depending only on $\Ii$, such that $b_{i,j}\in\Ii_0$ when $B_{i,j}\cdot F_i>0$, and $m_{i,j}\in\Ii_0$ when $\Mm_{i,X_i'}\cdot F_i>0$. 
By our assumption, if $b_{i,1}$ is strictly increasing, then $B_{i,1}$ is horizontal$/Z_i$, so $B_{i,1}\cdot F_i>0$, a contradiction; and if $b_{i,1}$ is not strictly increasing, then $m_{i,1}$ is strictly increasing and $\Mm_{i,1,X_i}$ is ample$/Z_i$, so $\Mm_{i,1,X_i}\cdot F_i>0$, a contradiction.

\medskip

\noindent\textbf{Step 2}. We prove the theorem when $(X_i,B_i,\Mm_i)$ is not klt for infinitely many $i$, and reduce to the case when $(X_i,B_i,\Mm_i)$ is $\Qq$-factorial klt and $\rho(X_i)=1$ for any $i$.

Since $\Mm_{i,j}\not\equiv\bm{0}$, $\Mm_{i,j,X_i}\not\equiv 0$ for any $i,j$. Let $D_i:=B_{i,1}$ if $b_{i,1}$ is strictly increasing and let $D_i:=\Mm_{i,1,X_i}$ if $b_{i,1}$ is not strictly increasing. Then $K_{X_i}+B_i+\Mm_i-\epsilon D_i$ is not pseudo-effective, so we may run a $(K_{X_i}+B_i+\Mm_i-\epsilon D_i)$-MMP for some $0<\epsilon\ll 1$, which terminates with a Mori fiber space $X_i'\rightarrow Z_i$ as $X_i$ is a surface. Let $B_i',B_{i,j}'$ and $D_i'$ be the images of $B_i,B_{i,j}$ and $D_i$ on $X_i'$ for any $i,j$, then $(X_i',B_i',\Mm_i)$ is $\Qq$-factorial lc, $X_i'$ is klt, and $K_{X_i'}+B_i'+\Mm_{i,X_i'}\equiv 0$. By \textbf{Step 1}, possibly passing to a subsequence, we may assume that $Z_i$ is a point of each $i$. Thus $\rho(X_i')=1$ for each $i$. Moreover, $(X_i',B_i'-\epsilon B'_{i,1},\Mm_i)$ is $\Qq$-factorial dlt if $b_{i,1}$ is strictly increasing, and $(X_i',B_i',\Mm_i-\epsilon\Mm_{i,1})$ is $\Qq$-factorial dlt if $b_{i,1}$ is not strictly increasing.

Assume that $\lfloor B_i'\rfloor\not=0$ for any $i$. Let $S_i$ be an irreducible component of $\lfloor B_i'\rfloor$, then $S_i$ is normal as  $(X_i',B_i'-\epsilon B'_{i,1},\Mm_i)$ or $(X_i',B_i',\Mm_i-\epsilon\Mm_{i,1})$ are $\Qq$-factorial dlt. If $b_{i,1}$ is strictly increasing, then since $\rho(X_i')=1$, $B'_{i,1}\cdot S_i>0$, and we get a contradiction to Lemma \ref{lem: adjunction} by considering adjunction to $S_i$. If $b_{i,1}$ is not strictly increasing, then $m_{i,1}$ is strictly increasing. Since $\rho(X_i')=1$ and $\Mm_{i,1}\not\equiv 0$, $\Mm_{i,1,X_i'}$ is ample, so $\Mm_{i,1,X_i'}\cdot S_i>0$, and we get a contradiction to Lemma \ref{lem: adjunction} again by considering adjunction to $S_i$. Thus, possibly passing to a subsequence, we may assume that $\lfloor B_i'\rfloor=0$ for any $i$.

Assume that $(X_i',B_i',\Mm_i)$ is not klt for any $i$. Let $E_i$ be an lc place of $(X_i',B_i',\Mm_i)$, then $E_i$ is exceptional$/X_i$. Let $f_i: \bar X_i\rightarrow X_i$ be the extraction of $E_i$, and $\bar B_i$, $\bar B_{i,j}$ the strict transforms of $B_i',B_{i,j}'$ on $\bar X_i$. Let $\bar D_i:=\bar B_{i,1}$ if $b_{i,1}$ is strictly increasing and let $\bar D_i:=\Mm_{i,1,\bar X_i}$ if $b_{i,1}$ is not strictly increasing. Since $\lfloor B_i'\rfloor\not=0$, $(X_i',B_i'-\epsilon B'_{i,1},\Mm_i)$ or $(X_i',B_i',\Mm_i-\epsilon\Mm_{i,1})$ is $\Qq$-factorial dlt, we have that $(X_i',B_i'-\epsilon B'_{i,1},\Mm_i)$ or $(X_i',B_i',\Mm_i-\epsilon\Mm_{i,1})$ is $\Qq$-factorial klt. Therefore, $\bar D_i|_{E_i}\not\equiv 0$. We get a contradiction to Lemma \ref{lem: adjunction} again by considering the adjunction to $E_i$.

Therefore, possibly passing to a subsequence, we may assume that $(X_i', B_i',\Mm_i)$ is klt for any $i$. Possibly replacing $(X_i,B_i,\Mm_i)$ with $(X_i',B_i',\Mm_i)$, we may assume that $(X_i,B_i,\Mm_i)$ is $\Qq$-factorial klt and $\rho(X_i)=1$ for any $i$.

\medskip

\noindent\textbf{Step 3}. Let $a_i:=\mld(X_i)$ be the minimal log discrepancy of $X_i$. In this step, we prove the case when $\lim_{i\rightarrow+\infty}a_i\not=0$.

Suppose that $\lim_{i\rightarrow+\infty}a_i\not=0$. Possibly passing to a subsequence, we may assume that $a:=\lim_{i\rightarrow+\infty}a_i$ is well-defined and $a>0$. By \cite[0.4(1)]{Ale94}, $X_i$ belongs to a bounded family. In particular, there exists a positive integer $N$, and very ample divisors $A_i$ on $X_i$, such that $-K_{X_i}\cdot A_i\leq N$ for any $i$. Possibly passing to a subsequence, we may assume that $M:=-K_{X_i}\cdot A_i$ is a constant positive integer. Therefore,
$$M=\sum b_{i,j}(B_{i,j}\cdot A_i)+\sum m_{i,j}(\Mm_{i,j,X_i}\cdot A_i).$$
Since $B_{i,j}\cdot A_i$ and $\Mm_{i,j,X_i}\cdot A_i$ are positive integers, $b_{i,j}$ and $m_{i,j}$ belong to a finite set. A contradiction.

\medskip

\noindent\textbf{Step 4}. In this step, we prove the case when there are two prime divisors $E_i, C_i$ on $X_i$, such that $\lim_{i\rightarrow+\infty}\mult_{E_i}B_i=\lim_{i\rightarrow+\infty}\mult_{C_i}B_i=1$.

We let $e_i:=\mult_{E_i}B_i$, $c_i:=\mult_{C_i}B_i$, and $\tilde B_i:=B_i-e_iE_i-c_iC_i$. Since $\rho(X_i)=1$, $E_i\equiv \mu_iC_i$ for some $\mu_i>0$. Possibly passing to a subsequence and switching $E_i, C_i$, we may assume that $\mu_i<1$ for each $i$. We let
$$\tilde B_i':=\tilde B_i+E_i+(1-c_i-e_i\mu_i)C_i.$$
By Lemma \ref{lem:surface ACC}, possibly passing to a subsequence, we have that the coefficients of $\tilde B_i$ belong to a DCC set, $(X_i,\tilde B_i',\Mm_i)$ is lc but not klt, $X_i$ is $\Qq$-factorial klt, and $K_{X_i}+\tilde B_i'+\Mm_{i, X_i}\equiv 0$. This contradicts \noindent\textbf{Step 2}.

\medskip

\noindent\textbf{Step 5}. Let $E_i$ be a prime divisor over $X_i$ such that $a(E_i,X_i,B_i,\Mm_i)=\mld(X_i,B_i,\Mm_i)$. In this step, we reduce to the case when $E_i$ is on $X_i$. Since $a(E_i, X_i, B_i,\Mm_i)<\mld(X_i)$, possibly passing to a subsequence, we may assume that $e_i$ is strictly decreasing, $e_i<1$ for each $i$, and $\lim_{i\rightarrow+\infty}e_i=0$. By  Lemma \ref{lem:extract divisor}, there exists an extraction $f_i: \bar X_i\rightarrow X_i$ of $E_i$. Let
$$K_{\bar X_i}+\bar B_i+e_iE_i+\Mm_{i,\bar X_i}:=f_i^*(K_{X_i}+B_i+\Mm_{i,X_i})$$
where $e_i:=1-a(E_i,X_i,B_i,\Mm_i)$, and let $0<e_i'<e_i$ be a real number such $\lim_{i\rightarrow+\infty}(e_i-e_i')=0$. Since 
$$K_{\bar X_i}+\bar B_i+e_iE_i+\Mm_{i,\bar X_i}\equiv 0,$$
$K_{\bar X_i}+\bar B_i+e_i'E_i+\Mm_{i,\bar X_i}$ is not pseudo-effective, so we may run a $(K_{\bar X_i}+\bar B_i+e_i'E_i+\Mm_{i,\bar X_i})$-MMP. Let $h_i: \bar X_i\rightarrow \bar X_i'$ be the first step of this MMP. Note that since $\bar X_i$ is a surface, $h_i$ is either a Mori fiber space or a divisorial contraction. Since $E_i$ is ample$/X_i'$, by \noindent\textbf{Step 1}, $h_i$ is not a Mori fiber space. Thus $h_i$ is a divisorial contraction of a prime divisor $F_i$. We have that $\rho(\bar X_i')=1$. Let $\bar B_i',E_i'$ be the image of $\bar B_i,E_i$ on $\bar X_i'$ for each $i$. 

By \textbf{Step 3}, possibly passing to a subsequence, there exists a prime divisor $C_i$ over $\bar X_i'$ such that $a(C_i,\bar X_i')$ is strictly decreasing and $\lim_{i\rightarrow+\infty}a(C_i,\bar X_i')=0$.  In particular, possibly passing to a subsequence, we may assume that $C_i$ is exceptional$/\bar X_i'$ and $a(C_i,\bar X_i',\bar B_i'+e_iE_i',\Mm_i)$ is strictly decreasing. Let $c_i:=1-a(C_i,\bar X_i',\bar B_i'+e_iE_i',\Mm_i)$ and let $g_i: \widehat X_i\rightarrow\bar X_i'$ be the extraction of $C_i$. We have
$$K_{\widehat X_i}+\widehat B_i+e_i'\widehat E_i+c_iC_i+\Mm_{i,\widehat X_i}=g_i^*(K_{\bar X_i'}+\bar B_i'+e_i'E_i'+\Mm_{i,\bar X_i'}),$$
where $\widehat B_i,\widehat E_i$ are the strict transforms of $\bar B_i',E_i'$ on $\widehat X_i$.

We run a $(-C_i)$-MMP and let $\phi_i:\widehat X_i\rightarrow\widehat X_i'$ be the first step of this MMP. Note that since $\widehat X_i$ is a surface, $\phi_i$ is either a Mori fiber space or a divisorial contraction. Since $C_i$ is ample$/\widehat X_i'$, by \noindent\textbf{Step 1}, $\phi_i$ is not a Mori fiber space, so it is a divisorial contraction. Let $\widehat B_i',\widehat E_i',C_i'$ be the images of $\widehat B_i,\widehat E_i,C_i$ on $\widehat X_i'$ respectively.

In particular, $\rho(\widehat X_i')=1$. By \noindent\textbf{Step 4}, $\widehat E_i$ is contracted by $\phi_i$. Then $(\widehat X_i',\widehat B_i'+c_iC_i',\Mm_i)$ is $\Qq$-factorial klt and 
$$K_{\widehat X_i'}+\widehat B_i'+c_iC_i'+\Mm_{i,X_i'}\equiv 0.$$
By Lemma \ref{lem:surface ACC}, possibly passing to a subsequence, $(\widehat X_i',\widehat B_i'+C_i',\Mm_i)$ is lc. Thus
$$K_{\widehat X_i}+\widehat B_i+e_i''\widehat E_i+C_i+\Mm_{i,\widehat X_i}=\phi_i^*(K_{\widehat X_i'}+\widehat B_i'+C_i'+\Mm_{i,X_i'})$$
for some $e_i''\leq 1$. Therefore, 
$$(K_{\widehat X_i}+\widehat B_i+\widehat E_i+C_i+\Mm_{i,\widehat X_i})\cdot \widehat E_i\leq 0.$$
By  Lemma \ref{lem:surface ACC} again, possibly passing to a subsequence, $(\bar X_i',\bar B_i'+E_i',\Mm_i)$ is lc for any $i$. Thus
$$K_{\widehat X_i}+\widehat B_i+\widehat E_i+c_i'C_i+\Mm_{i,\widehat X_i}=g_i^*(K_{\bar X_i'}+\bar B_i'+E_i'+\Mm_{i,\bar X_i'})$$
for some $c_i'\leq 1$. Therefore, 
$$(K_{\widehat X_i}+\widehat B_i+\widehat E_i+C_i+\Mm_{i,\widehat X_i})\cdot C_i\leq 0.$$
Since $\rho(\widehat X_i)=2$ and $C_i,\widehat{E}_i$ are negative extremal rays in $\overline{NE}(X_i)$, 
$$-(K_{\widehat X_i}+\widehat B_i+\widehat E_i+C_i+\Mm_{i,\widehat X_i})$$ is nef. Thus
$$(1-c_i)C_i'\equiv (\phi_i)_*(-(K_{\widehat X_i}+\widehat B_i+\widehat E_i+C_i+\Mm_{i,\widehat X_i}))$$
is pseudo-effective. This is not possible.
\end{proof}

\subsection{ACC on threefolds}

\begin{prop}\label{prop:ACC on threefolds}
    Conjecture \ref{conj: ACC for lc thresholds} holds when $d = 3$.
\end{prop}
\begin{proof}
    We may assume that $1\in\Ii$. Suppose the proposition is false, then there exists a sequence $(X_i, B_i,\Mm_i; D_i,\Nn_i)$ such that $\dim X_i=3$, $B_i, D_i\in\Ii$, $\Mm,\Nn\in\Nef(X,\Ii)$, and $t_i:=\lct(X_i, B_i,\Mm_i; D_i,\Nn_i)$ is strictly increasing. In particular, we may assume that $t_i>0$ for each $i$. 

    Suppose that for infinitely many $i$, there exists an lc place $E_i$ of $(X_i,\Delta_i:=B_i+t_iD_i,\Pp_i:=\Mm_i+t_i\Nn_i)$ on $X_i$ such that
    $$a(E_i,X_i,B_i+(t_i+\epsilon)D_i,\Mm_i+(t_i+\epsilon)\Nn_i)<0$$
    for any $\epsilon>0$. Then $E_i$ is a component of $B_i+t_iD_i$, $\mult_{E_i}D_i>0$, and $\mult_{E_i}(B_i+t_iD_i)=1$, which is not possible as $\Ii$ satisfies the DCC and $t_i$ is strictly increasing. Thus, possibly passing to a subsequence, we may assume that no such $E_i$ exists. 
    
    By Lemma \ref{lem: special extraction}, there are two projective birational morphisms $h_i: X_i'\rightarrow X_i$ and $g_i: Y_i\rightarrow X_i'$, such that
    \begin{itemize}
        \item   $h_i$ is a $\Qq$-factorial dlt modification of $(X_i,\Delta_i,\Pp_i)$,
        \item any $h_i$-exceptional prime divisor $P_i$ is an lc place of $(X_i,B_i,\Mm_i)$,
        \item $g_i$ extracts a unique prime divisor $E_i$ such that $-E_i$ is ample$/X_i'$, 
        \item $a(E_i,X_i,B_i+t_iD_i,\Mm_i+t_i\Nn_i)=0$ and $a(E_i,X_i,B_i,\Mm_i)>0$, and
        \item $(Y_i,B_{Y_i},\Mm)$ is $\Qq$-factorial dlt, where $B_{Y_i}$ is the strict transform of $B_i$ on $Y_i$ plus the reduced $(h_i\circ g_i)$-exceptional divisor. By Lemma \ref{lem: dlt implies normal}, $E_i$ is normal.
    \end{itemize}
    Let $D_i',D_{Y_i}$ be the strict transforms of $D_i$ on $X_i'$, $Y_i$ respectively, and let $B_i'$ be the strict transform of $B_i$ on $X_i'$ plus the reduced $h_i$-exceptional divisor. By our construction,
    $$t_i=\lct(X_i',B_i',\Mm_i;D_i',\Nn_i),$$
    and $D_{Y_i}+\Nn_{i,Y_i}=g_i^*(D_i'+\Nn_{i,X_i'})+b_iE_i$ for some $b_i>0$. Therefore, $(D_{Y_i}+\Nn_{i,Y_i})|_{E_i}$ is ample. 

    Let $x_i':=\Center_{X_i'}E_i$, and let $F_i$ be the reduced variety associated to a general fiber of $E_i\rightarrow x_i'$,
    $$K_{E_i}+B_{E_i}(t)+\Mm_i(t)_{E_i}:=(K_{Y_i}+B_{Y_i}+tD_{Y_i}+\Mm_{i,Y_i}+t\Nn_{i,Y_i})|_{E_i}$$
    for any real number $t$. 
    When $\dim F_i = 2$, then $K_{F_i} = K_{E_i}$ and $$K_{F_i}+B_{F_i}(t_i)+\Mm^{F_i}(t_i)_{F_i}\equiv 0$$
    and
    $$K_{F_i}+B_{F_i}(0)+\Mm^{F_i}(0)_{F_i}$$
    is anti-ample. By Lemma \ref{lem:global ACC}, the coefficients of $B_{E_i}(t_i)\in \Gamma_0$ and $\Mm_i(t_i) \in \Nef^0(\Gamma_0)$ for some finite set $\Gamma_0$ depending only on $\Ii$.
    When $\dim F_i = 1$, then 
    $$(K_{E_i} + F_i) \cdot F_i = 2p_a(F_i)-2, F_i^2=0.$$ 
    Writing $B_{E_i}(t_i) = \sum_j b_{i,j} B_{i,j}, M_i(t_i)_{E_i} = \sum_j m_{i,j} M_{i,j}$, we have 
    $$0 = (K_{E_i} +B_{E_i}(t_i) + M_i(t_i)_{E_i}) \cdot F_i = 2p_a(F_i)-2 + \sum_j b_{i,j}(B_{i,j} \cdot F_i) + \sum_j m_{i,j} (M_{i,j} \cdot F_i)$$ 
    and $B_{i,j} \cdot F_i, M_{i,j} \cdot F_i \in \mathbb{N}$. 
    Then there are also only finitely many choices of $b_{i,j}, m_{i,j}$ as they belong to a DCC set.
    After expanding $\Gamma_0$, we may assume that the coefficients of $B_{E_i}(t_i)\in \Gamma_0$ and $\Mm_i(t_i) \in \Nef^0(\Gamma_0)$ for all $i$.
    By Lemma \ref{lem: adjunction}(1), $K_{F_i}+B_{F_i}(t)+\Mm^{F_i}(t)_{F_i}$ is constant with respect to $t$ for some $i$.
    Then $K_{F_i}+B_{F_i}(t)+\Mm^{F_i}(t)_{F_i} \equiv 0$ for all $t \in \Rr$,
    which is impossible as $K_{F_i}+B_{F_i}(0)+\Mm^{F_i}(0)_{F_i}$ is anti-ample.
\end{proof}

\begin{rem}
    We essentially get the above results by computing the intersection number of $F_i$ with the divisor, which works in the surface case, not by applying the adjunction formula. It would be worthwhile to investigate the adjunction formula for general fibers of fibrations, or Mori fiber spaces, in characteristic $p$. We believe that one can get more general results following the ideas in \cite{PW22}.
\end{rem}

\section{Existence of Shokurov-type polytope}

\begin{prop}\label{prop:length of extremal rays}
    Let $(X, B,\Mm)/U$ be a $\Qq$-factorial lc generalized pair of dimension $\leq 3$ such that $X$ is klt. Then any $(K_X+B+\Mm_X)$-negative extremal ray $R/U$ is spanned by a rational curve $C$ such that 
    $$0<-(K_X+B+\Mm_X)\cdot C\leq 6.$$
\end{prop}

\begin{proof}
Let $R$ be a $(K_X+B+\Mm_X)$-negative extremal ray$/U$. Then $R$ is a $(K_X+B+\Mm_X+A)$-negative extremal ray for some ample$/U$ $\Rr$-divisor $A$. By Lemma \ref{lem: hl22 3.4 char p}, $R$ is a $(K_X+\Delta)$-negative extremal ray$/U$ for some klt pair $(X,\Delta)$ such that $K_X+\Delta\sim_{\mathbb R,U}K_X+B+\Mm_X+A$. The proposition follows from \cite[Theorem 1.3(4)]{HNT20}.
\end{proof}

\begin{proof}[Proof of Theorem \ref{thm: shokurov polytope intro}] Let $B(\bm{v}):=\sum_{i=1}^m v_iB_i$ and $\Mm(\bm{v}):=\sum_{i=m+1}^{m+n}v_i\Mm_i$ for any $\bm{v}:=(v_1,\dots,v_{m+n})\in\mathbb R^{m+n}$. By the existence of log resolutions, there exists an open subset $U_1\ni\bm{v}_0$ in the rational envelope of $\bm{v}_0$, such that for any $\bm{v}\in U_1$, $(X, B(\bm{v}),\Mm(\bm{v}))$ is lc. We let $c:=\dim U_1$ and let $\bm{v}_1,\dots,\bm{v}_{c+1}$ be vectors in $U_1\cap\mathbb Q^{m+n}$ such that $\bm{v}_0$ is contained in the convex hull $U_2$ spanned by
    $\bm{v}_1,\dots,\bm{v}_{c+1}$. Then there exist positive real numbers $a_1,\dots,a_{c+1}$ such that $\sum_{i=1}^{c+1}a_i\bm{v}_i=\bm{v}_0$ and $\sum_{i=1}^{c+1}a_i=1$. We let $I$ be a positive integer such that $I(K_{X}+B(\bm{v}_i)+\Mm(\bm{v}_i)_X)$ is Cartier for each $i$. Let $a_0:=\min_{1\leq i\leq c+1}\{a_i\}$. 

    Consider the set
    $$\Ii:=\left\{\sum a_i\gamma_i\mid\gamma_i\in [-6I,+\infty)\cap\mathbb Z\right\}\cap (0,+\infty).$$
    We have $\gamma_0:=\inf\{\gamma\in\Ii\}>0$. We let $U_0$ be the interior of the set
    $$\left\{\frac{1}{6I+\gamma_0}(6I\bm{v}_0+\gamma_0\bm{v})\Bigg| \bm{v}\in U_2\right\}.$$

    We show that $U_0$ satisfies our requirement. By our construction, $(X,B(\bm{v}),\Mm(\bm{v}))$ is lc for any $\bm{v}\in U_0$, so we only need to show that $K_{X}+\sum_{i=1}^mv_iB_i+\sum_{i=m+1}^{m+n}v_i\Mm_{i,X}$ is nef$/U$ for any $\bm{v}=(v_1,\dots,v_{m+n})\in U_0$. We let $R$ be an extremal ray in $\overline{NE}(X/U)$. There are three cases.

    \medskip

    \noindent\textbf{Case 1}. $(K_{X}+B+\Mm_X)\cdot R=0$. In this case, by \cite[Lemma 5.3]{HLS19}, $(K_{X}+B(\bm{v})+\Mm(\bm{v})_X)\cdot R=0$ for any $\bm{v}\in U_1$, so $(K_{X}+B(\bm{v})+\Mm(\bm{v})_X)\cdot R=0$ for any $\bm{v}\in U_0$.

    \medskip

    \noindent\textbf{Case 2}. $(K_{X}+B(\bm{v}_i)+\Mm(\bm{v}_i)_X)\cdot R\geq 0$ for any $i$. In this case, $(K_{X}+B(\bm{v})+\Mm(\bm{v})_X)\cdot R\geq 0$ for any $\bm{v}\in U_2$, so $(K_{X}+B(\bm{v})+\Mm(\bm{v})_X)\cdot R\geq 0$ for any $\bm{v}\in U_0$.

    \medskip

    \noindent\textbf{Case 3}. $(K_{X}+B+\Mm_X)\cdot R>0$ and $(K_{X}+B(\bm{v}_j)+\Mm(\bm{v}_j))\cdot R<0$ for some $j$. In this case, by Proposition \ref{prop:length of extremal rays},  $R$ is spanned by a curve $C$ such that $(K_{X}+B(\bm{v}_i)+\Mm(\bm{v}_i)_X)\cdot C\geq -6$ for any $i$. Thus $$I(K_{X}+B(\bm{v}_i)+\Mm(\bm{v}_i)_X)\cdot C\in [-6I,+\infty)\cap\mathbb Z,$$
    so
    $$I(K_{X}+B+\Mm_X)\cdot C\in\Ii.$$
    Then for any $\bm{v}\in U_0$, there exists $\bm{v}'\in U_2$ such that $(6I+\gamma_0)\bm{v}=6I\bm{v}_0+\gamma_0\bm{v}'$. We have
    \begin{align*}
        &I(K_{X}+B(\bm{v})+\Mm(\bm{v})_X)\cdot C\\
        =&\frac{\gamma_0}{6I+\gamma_0}I(K_{X}+B(\bm{v}')+\Mm(\bm{v}')_X)\cdot C+\frac{6I}{6I+\gamma_0}I(K_{X}+B(\bm{v}_0)+\Mm(\bm{v}_0)_X)\cdot C\\
        \geq &\frac{\gamma_0}{6I+\gamma_0}\cdot (-6I)+\frac{6I}{6I+\gamma_0}\cdot\gamma_0=0,
    \end{align*}
    so $I(K_{X}+B(\bm{v})+\Mm(\bm{v})_X)\cdot R\geq 0$. The theorem follows.
\end{proof}

\section{Termination of flips in terminal case}


\begin{defn}[Difficulty]
    Let $(X,B,\Mm)/U$ be a $\Qq$-factorial NQC terminal generalized pair with $B = \sum b_i B_i$ and $\Mm = \sum \mu_j \Mm_j$, where $B_i$ are the irreducible components of $B$ and each $\Mm_j$ is a nef$/U$ $\bb$-Cartier $\bb$-divisor. Let $b:=\max\{b_i,0\}$ and let
    \[ \Ss_b := \left\{ \sum n_ib_i + \sum m_j \mu_j \big| n_i,m_j \in \Zz_{\geq 0} \right\} \cap [b,+\infty).\]
For any $\xi\in\Ss_b$, we define
\[d_\xi(X,B,\Mm) := \#\{E\mid E \text{ is exceptional}/X, a(E,X,B,\Mm) < 2 - \xi\}.\]
    The \emph{difficulty} of $(X,B,\Mm)$ is defined as 
    \[ d(X,B,\Mm) := \sum_{\xi \in \Ss_b} d_\xi(X,B,\Mm).\]
\end{defn}

\begin{lem}\label{lem:finiteness of difficulty}
Let $(X,B,\Mm)/U$ be a $\Qq$-factorial NQC terminal generalized pair. Then
    \[ 0 \leq d(X,B,\Mm) < +\infty.\]
\end{lem}
\begin{proof}
    We follow the proof of \cite[Lemma 2.22]{CT23}. It still holds in the positive characteristic. For convenience, we provide a full proof.
    
    For $\xi\geq 1$, $d_\xi(X,B,\Mm) = 0$ since $(X,B,\Mm)$ is terminal.
    Since the set $\Ss_b\cap [0,1)$ is finite, we only need to show that $d_\xi(X,B,\Mm)<\infty$ for any fixed $\xi$. Since $d_\xi(X,B,\Mm)\leq d_b(X,B,\Mm)$ for any $\xi \in [0,1)$, we only need to show that $d_b(X,B,\Mm) < \infty$.
    Let $f: X' \to X$ be a log resolution of $(X,B)$ such that $\Mm$ descends to $X'$, and let 
    $$ K_{X'}+B'+\Mm_{X'}:=f^*(K_X+B+\Mm_X).$$
    We may write $B' = (B')^+ - (B')^-$, where $(B')^+\geq 0,(B')^-\geq 0$, and $(B')^+\wedge (B')^-=0$. Possibly after further blow-ups, we may assume that $\Supp (B')^+$ is non-singular. Since $(X,B,\Mm)$ is terminal, $(B')^+=f^{-1}_*B$. Thus for any prime divisor $F$ that is exceptional$/X$, we have $$a(F,X,B,\Mm) = a(F,X',B',\Mm')=a(F,X',B') \geq 2 - b $$ by \cite[Corollary 2.31(3)]{KM98}. Therefore, 
    $$d_b(X,B,\Mm)\leq\#\{E\mid E\text{ is exceptional}/X, E\text{ is on }X'\}<+\infty.$$
\end{proof}

\begin{lem}\label{lem:blow-up disc}
Let $(X,B,\Mm)/U$ be an NQC generalized pair such that $B=\sum_i b_iB_i$ and $\Mm=\sum_j\mu_j\Mm_j$, where $B_i$ are the irreducible components of $B$, each $\Mm_j$ is a nef$/U$ $\bb$-Cartier $\bb$-divisor, and each $\mu_j\geq 0$. Assume that $X$ is smooth at the generic point of a codimension $k \geq 2$ closed subvariety $V$ of $X$ and let $E$ be the irreducible component of the blow-up of $V$ which dominates $V$. Then
$$a(E,X,B,\Mm) = k - \sum_i n_i b_i - \sum_j m_j \mu_j$$ for some non-negative integers $n_i$ and $m_j$.
\end{lem}
\begin{proof}
Let $f: X'\rightarrow X$ be the blow-up of $V$. Then
$$a(E,X,B,\Mm)=a(E,X,B)-\mult_E(f^*\Mm_X-\Mm_{X'}).$$
By \cite[Lemma 2.29]{KM98}, $a(E,X,B) = k - \sum_i n_i b_i$ for some non-negative integers $n_i$. By the negativity lemma, $f^*\Mm_{j,X}- \Mm_{j,X'}\geq 0$ for each $j$. Since $\Mm_{j,X}$ is Cartier near the generic point of $V$, $m_j:=\mult_E(f^*\Mm_{j,X}- \Mm_{j,X'})$ is a non-negative integer and the lemma follows. \end{proof}

\begin{thm}\label{thm: termination of terminal case}
    Let $(X,B,\Mm)/U$ be a 3-dimensional $\Qq$-factorial NQC terminal generalized pair.
    Then any sequence of $(K_X+B+\Mm_X)$-flips$/U$ terminates.
\end{thm}
\begin{proof}
 Let
    $$(X,B,\Mm)=:(X_1,B_1,\Mm)\dashrightarrow (X_2,B_2,\Mm)\dashrightarrow (X_3,B_3,\Mm)\dashrightarrow\dots\dashrightarrow (X_i,B_i,\Mm)\dashrightarrow\dots$$
    be a sequence of $(K_X+B+\Mm_X)$-flips$/U$ and let $f_i: X_i\rightarrow W_i$ be the flipping contractions. Since log discrepancies do not decrease under the MMP and $X\dashrightarrow X_i$ does not contract any divisor for any $i$, $(X_i, B_i,\Mm)$ is terminal for each $i$. We prove by induction on the number of components of $B_1$.

\medskip

\noindent\textbf{Step 1.}  First suppose $B_1 = 0$. Since $X_{i+1}$ is terminal, it is smooth along the generic point of a flipped curve $C_{i+1}$. Let $\eta_{i+1}$ be the generic point of $C_{i+1}$. Let $E_{i+1}$ be the exceptional divisor obtained by blowing up $X_{i+1}$ at $\eta_{i+1}$, then by Lemma \ref{lem:blow-up disc}, 
$$a(E_{i+1}, X_{i}, \Mm_{i}) < a(E_{i+1}, X_{i+1}, \Mm_{i+1})=2 - \sum_j m_j \mu_j$$ 
for some non-negative integers $m_j$. Therefore, $d(X_i,M_i) > d(X_{i+1},M_{i+1})$. Since $d(X_1, M_1) < + \infty$, the sequence of flips cannot be infinite.

\medskip

\noindent\textbf{Step 2.}  Now suppose that any sequence of flips terminates when the number of components of $B_1 \leq s - 1$ for some positive integer $s$. If the number of components of $B_1 = s$, then we may write $B_1=\sum_{k=1}^sb_kB_{1,k}$ where $B_{1,k}$ are the irreducible components of $B_1$, and let $b:=\max_kb_k$. We define
$$D_1:=\sum_{b_k=b}B_{1,k}$$
and let $D_i$ be the image of $D_1$ on $X_i$ for each $i$.

First, we reduce to the case when $D_i$ does not contain any flipped curve of $X_i\rightarrow W_{i-1}$. Suppose that $D_i$ contains a flipped curve $C_i$. Since $X_{i}$ is terminal, it is smooth along the generic point of a flipped curve $C_{i}$. Let $\eta_{i}$ be the generic point of $C_{i}$. Let $E_{i}$ be the exceptional divisor obtained by blowing up $X_{i}$ at $\eta_{i}$, then by Lemma \ref{lem:blow-up disc}, then 
$$a(E_i,X_{i-1},B_{i-1},\Mm)<a(E_i,X_i,B_i,\Mm)=2-\sum n_ib_i-\sum m_j\mu_j$$
for some non-negative integer $n_i,m_j$, and when $b_k=b$, we have $n_k>0$. Therefore, 
$$2-\sum n_ib_i-\sum m_j\mu_j\leq 2-\sum_{b_k=b}n_kb_k\leq 2-b.$$
So $d(X_{i+1},B_{i+1},\Mm)>d(X_i,B_i,\Mm)$. Thus, possibly after truncating, we may assume that no $D_i$ contains any flipped curve of $X_i\rightarrow W_{i-1}$.

Next, we reduce to the case when $D_i$ does not contain any flipping curve of $X_i\rightarrow W_{i}$. Let $\nu_i: D_i^\nu\rightarrow D_i$ be the normalization of $D_i$ and $\varphi_i: D_i^\nu\rightarrow W_i^\nu$, $\varphi_i^+: D_{i+1}^\nu\rightarrow W_i^\nu$ the induced contractions. Let $g_i: D_i^\nu\dashrightarrow D_{i+1}^\nu$ be the induced birational map, then since no flipped curve of $X_{i+1}\rightarrow W_i$ is contained in $D_{i+1}$, $g_i$ does not extract any divisor, so $g_i$ is a contraction. We let $q_1: D_1^{\mu}\rightarrow D_1^{\nu}$ be the minimal resolution of $D_1^{\nu}$ and let $q_i: D_1^{\mu}\rightarrow D_i^{\nu}$ be the induced contraction and let
$$\rho_i:=\rho(D_1^{\mu}/U)-\#\{E\mid E\text{ is a }q_i\text{-exceptional prime divisor}\}.$$
Then $\rho_i$ is non-increasing, and $\rho_i\geq 0$ by \cite[Lemma 1.6]{AHK07} (note that \cite{AHK07} works over $\mathbb C$, but the same lines of the proof of \cite[Lemma 1.6]{AHK07} work over an algebraically closed field of any characteristic). Thus, possibly truncating the MMP, we may assume that $\rho_i$ is a constant. Thus $g_i$ does not contract any divisor. However, any flipping curve of $X_i\rightarrow W_i$ that is contained in $D_i$ is contracted by $g_i$. Therefore, $D_i$ does not contain any flipping curve of $X_i\rightarrow W_i$.

Hence, $D_i \cdot C \geq 0$ for every flipping curve $C \subseteq X_i$. Therefore, any $(K_{X_i} + B_i + \Mm_{i,X_i})$-flip is also a $(K_{X_i} + (B_i - b D_i) + \Mm_{i,X_i})$-flip, but $B_i-b D_i$ contains at most $s-1$ irreducible components, and we are done by induction on $s$.
\end{proof}

\section{Existence of log minimal model}

\subsection{Log minimal model}

\begin{lem}\label{lem: log minimal model for log smooth model}
    Let $(X,B,\Mm)/U$ be an lc generalized pair of dimension $\leq 3$, $f: W\rightarrow X$ a log resolution of $(X,B)$ such that $\Mm$ descends to $W$, and $B_W:=f^{-1}_*B+\Exc(f)$. Let $(X',B',\Mm)/U$ be a log minimal model of $(W,B_W,\Mm)/U$. Then $(X',B',\Mm)/U$ is also a log minimal model of $(X,B,\Mm)/U$.
\end{lem}
\begin{proof}
    By definition, we only need to show that 
    \begin{enumerate}
        \item for any prime divisor $D$ on $X$ that is exceptional$/X'$, we have $a(D,X,B,\Mm)<a(D,X',B',\Mm)$, and
        \item $(X',B',\Mm)/U$ is a log birational model of $(X,B,\Mm)/U$.
    \end{enumerate}
First, we prove (1). For any prime divisor $D$ on $X$ that is exceptional$/X'$, $\Center_WD$ is a divisor, so $a(D,W,B_W,\Mm)<a(D,X',B',\Mm)$. Since $(X,B,\Mm)$ is lc, $a(D,X,B,\Mm)\leq a(D,W,B_W,\Mm)$. This implies (1).

Now we prove (2). For any prime divisor $D$ on $X'$ that is exceptional$/X$, if $D$ is exceptional$/W$, then since $(X',B',\Mm)/U$ is a log birational model of $(W,B_W,\Mm)/U$, $\mult_DB'=1$. If $D$ is not exceptional$/W$, then $D_W:=\Center_WD$ is a divisor and $D$ is exceptional$/X$, so $\mult_DB'=\mult_{D_W}B_W=1$. This implies (2).
\end{proof}

\begin{lem}\label{lem: bir12 2.6}
    Let $(X,B,\Mm)/U$ be an lc generalized pair of dimension $\leq 3$, $(X',B',\Mm)/U$ a bs-weak lc model of $(X,B,\Mm)/U$, and $g: W\rightarrow X$ and $h: W\rightarrow X'$ two birational morphisms. Then
    $$g^*(K_X+B+\Mm_X)=h^*(K_{X'}+B'+\Mm_{X'})+F$$
    for some $F\geq 0$ that is exceptional$/X'$.
\end{lem}
\begin{proof}
    By the definition of bs-weak lc models, $g_*F\geq 0$. By the negativity lemma, $F\geq 0$. If $F$ is not exceptional$/X'$, then let $D$ be a component of $F$ that is not exceptional$/X$. If $D$ is not exceptional$/X$, then $a(D,X,B,\Mm)=a(D,X',B',\Mm)$ so $D$ is not a component of $F$, a contradiction. If $D$ is exceptional$/X$, then $\mult_DB'=1$. Since $F\geq 0$, 
    $$0\leq a(D,X,B,\Mm)\leq a(D,X',B',\Mm)=0,$$
    so $D$ is not a component of $F$, a contradiction.
\end{proof}

\begin{lem}\label{lem: bs wlcm implies lm}
    Let $(X, B,\Mm)/U$ be an lc generalized pair of dimension $\leq 3$. Assume that $(X, B,\Mm)/U$ admits a bs-weak lc model. Then $(X, B,\Mm)/U$ admits a log minimal model.
\end{lem}
\begin{proof}
    Let $(X',B',\Mm)/U$ be a bs-weak lc model of $(X,B,\Mm)/U$. Let $g: W\rightarrow X$ and $h: W\rightarrow X'$ be a common log resolution such that $\Mm$ descends to $W$. Let $B_W:=g^{-1}_*B+\Exc(g)$. Then
    $$K_W+B_W+\Mm_W=g^*(K_X+B+\Mm_X)+F_1$$
    for some $F_1\geq 0$ that is exceptional$/X$. By Lemma \ref{lem: bir12 2.6},
    $$g^*(K_X+B+\Mm_X)=h^*(K_{X'}+B'+\Mm_{X'})+F_2$$
    for some $F_2\geq 0$ that is exceptional$/X'$. Thus
    $$K_W+B_W+\Mm_W=h^*(K_{X'}+B'+\Mm_{X'})+F_1+F_2$$
    Since $F_1$ is exceptional$/X$, for any component $D$ of $F_1$ that is not exceptional$/X'$, $\mult_DB'=1$, so $\mult_DB_W>1$, which is not possible. Thus $F_1$ is exceptional$/X'$, so $F_1+F_2$ is exceptional$/X'$.
    We let $B_W'\geq B_W$ be the unique $\Rr$-divisor such that $\Supp(B_W'-B_W)\subset\Supp(F_1+F_2)$, and for any irreducible component $D$ of $F_1+F_2$, $\mult_DB_W'=1$. Then $(W,B_W')$ is log smooth and $$K_W+B_W'+\Mm_W\sim_{\mathbb R,X'}F\geq 0$$
    where $\Supp F=\Supp(F_1+F_2)$. We run a $(K_W+B_W'+\Mm_W)$-MMP$/X'$, which is an $F$-MMP. By Proposition \ref{prop: special termination}, the MMP terminates near $\lfloor B_W'\rfloor$. In particular, this MMP terminates near $\Supp F$, so this MMP terminates with a log minimal model $(X'',B'',\Mm)/X'$ of $(W,B_W',\Mm)/X'$. Let $F''$ be the image of $F$ on $X''$, then $K_{X''}+B''+\Mm_{X''}\sim_{\mathbb R,X'}F''\geq 0$, so $F''$ is nef$/X'$. By the negativity lemma, $F''=0$.    Let $f: X''\rightarrow X'$ be the induced birational morphism, then $K_{X''}+B''+\Mm_{X''}=f^*(K_{X'}+B'+\Mm_{X'})$. Therefore, $K_{X''}+B''+\Mm_{X''}$ is nef$/U$, so $(X'',B'',\Mm)/U$ is a log minimal model of $(W,B_W',\Mm)$.

    Since $F''=0$, any irreducible component of $B_W'-B_W$ is exceptional$/X''$. Since $B_W'\geq B_W$, $(X'',B'',\Mm)/U$ is a log minimal model of $(W,B_W',\Mm)$. By Lemma \ref{lem: log minimal model for log smooth model}, $(X'', B'',\Mm)/U$ is a log minimal model of $(X, B,\Mm)/U$.
 \end{proof}

\subsection{Weak Zariski decomposition}

\begin{defn}[Weak Zariski decomposition]
Let $X\rightarrow U$ be a projective morphism from a normal quasi-projective variety to a quasi-projective variety and let $D$ be an $\Rr$-Cartier $\Rr$-divisor on $X$. We say that $D$ admits a \emph{weak Zariski decomposition}$/U$ (resp. an \emph{NQC weak Zariski decomposition}$/U$) if there exists a projective birational morphism $f: W\rightarrow X$ such that
$$f^*D=P+N,$$
where $P$ is nef$/U$ (resp. NQC$/U$) and $N\geq 0$. $(f,P,N)$ is called a \emph{weak Zariski decomposition}$/U$ (resp. an \emph{NQC weak Zariski decomposition}$/U$) of $(X,B,\Mm)/U$.
\end{defn}

\begin{deflem}[MMP using weak Zariski decomposition]\label{deflem: MMP-Zariski}
Let $(X,B,\Mm)/U$ be a projective $\Qq$-factorial lc generalized pair of dimension $3$ such that $X$ is klt,
$$K_X+B+\Mm_X=P+N$$
for some nef$/U$ $\Rr$-divisor $P$ and $N\geq 0$, and $\Supp N\subset\lfloor B\rfloor$. Let
$$\mu:=\sup\{ t \in [0,1] | P + tN \text{ is nef}/U \}.$$ 
Then either $\mu=1$, or there exists an extremal ray$/U$ $R$ such that $(K_X + B + \Mm_X) \cdot R < 0$ and $(P + \mu N) \cdot R = 0$.

Moreover, let $\phi: (X,B,\Mm)\dashrightarrow (X',B',\Mm)$ be a divisorial contraction$/U$ or flip$/U$, such that it is a step of a $(K_X+B+\Mm_X)$-MMP$/U$ which contracts $R$. Then $(X',B',\Mm)/U$ is lc, $X'$ is $\Qq$-factorial klt, and we may write
$$K_{X'}+B'+\Mm_{X'}=P'+N'$$
for some nef$/U$ $\Rr$-divisor $P'$ and $N'\geq 0$, and $\Supp N'\subset\Supp\lfloor B'\rfloor$. We may replace $(X,B,\Mm)$ and $P,N$ with $(X',B',\Mm)$ and $P',N'$ respectively and repeat this process. Such a process will be called a \emph{$(K_X+B+\Mm_X)$-MMP$/U$ using weak Zariski decomposition $(P,N)$}.
\end{deflem}

\begin{proof}
First, we prove the existence of such $R$. We may assume that $\mu<1$. Replacing $P$ with $P + \mu N$ we may assume $\mu = 0$ and $N\not=0$. Then for every $\epsilon' > 0$, $P + \epsilon' N$ is not nef$/U$. In particular, for every $\epsilon' > 0$, there exists a $(K_X+B+\Mm_X)$-negative extremal ray$/U$ $R$ such that $(P + \epsilon' N) \cdot R < 0$, but $(P + \epsilon N) \cdot R = 0$ for some $0\leq\epsilon<\epsilon'$. 

Suppose that the lemma does not hold. Then we can find a strictly decreasing sequence $\{ \epsilon_i \}_{i=1}^{+\infty}$ such that $\lim_{i \to+\infty} \epsilon_i = 0$ and $(P + \epsilon_i N) \cdot R_i = 0$ for some $(K_X+B+\Mm_X)$-negative extremal ray$/U$ $R_i$. By Proposition \ref{prop:length of extremal rays}, each $R_i$ is generated by a curve $C_i$ such that 
$$0 < -(K_X+B+\Mm_X) \cdot C_i \leq 6.$$
We may assume that $\epsilon_1<1$. Since $\Supp N \subset \lfloor B \rfloor$, there exists $\delta > 0$ such that $(K_X+B+\Mm_X-\delta N) \cdot C_i< 0$ for each $i$, $B - \delta N \geq 0$, and $\Supp(B - \delta N) = \Supp B$. We have $$K_X+B+\Mm_X - \delta N= P+ (1-\delta) N.$$ 

We reduce to the case when $N\cdot C_i$ is a constant for any $i$. Let $\tau$ be a positive real number such that for any irreducible component $S$ of $N$, 
$$B\geq B-\delta N+\tau S\geq B-\delta N-\tau S\geq 0.$$
By Proposition \ref{prop:length of extremal rays}, for any $i$ and any irreducible component $S$ of $N$, if $S\cdot R_i\geq 0$, then 
$$0<-(K_X+B+\Mm_X-\delta N)\cdot C_i\leq -(K_X+B+\Mm_X-\delta N-\tau S)\cdot C_i\leq 6$$
so $0\leq S\cdot C_i<\frac{6}{\tau}$, and if $S\cdot R_i<0$, then
$$0<-(K_X+B+\Mm_X-\delta N)\cdot C_i<-(K_X+B+\Mm_X-\delta N+\tau S)\cdot C_i\leq 6,$$
so $0\leq -S\cdot C_i<\frac{6}{\tau}$. Let $I$ be a positive integer such that $IS$ is Cartier for any irreducible component $S$ of $N$, then $S\cdot C_i\in (-\frac{6}{\tau},\frac{6}{\tau})\cap\frac{1}{I}\mathbb Z$. Therefore, passing to a subsequence, we may assume that $S\cdot C_i$ is a constant for any irreducible component $S$ of $N$, hence $N\cdot C_i$ is a constant. Therefore,
$$\lim_{i\rightarrow+\infty}P\cdot C_i=\lim_{i\rightarrow+\infty}-\epsilon_iN\cdot C_i=0,$$
so possibly passing to a subsequence, $P\cdot C_i$ is strictly decreasing, hence
$$(K_X+B+\Mm_X)\cdot C_i=(P+N)\cdot C_i$$
is strictly decreasing.

There exist positive real numbers $a_1,\dots,a_k$ and lc generalized pairs $(X,B_j,\Mm_j)$, such that each $B_j$ is a $\Qq$-divisor, each $\Mm_j$ is a $\Qq$-$\bb$-divisor, $\sum_{j=1}^ka_j=1$, $\sum_{j=1}^ka_jB_ijB$, and $\sum_{j=1}^ka_j\Mm_j=\Mm$. Let $I$ be a positive integer such that $I(K_{X_j}+B_j+\Mm_{j,X})$ is Cartier for each $j$. Then for any $i,j$, by Proposition \ref{prop:length of extremal rays},
$$(K_{X}+B_j+\Mm_{j,X})\cdot C_i\in\frac{1}{I}\mathbb Z\cap [-6,+\infty),$$
so
$$(K_X+B+\Mm_X)\cdot C_i\in\Ii:=\left\{\sum_{j=1}^ka_j\gamma_j\Big| I\gamma_j\in [-6,+\infty)\right\}\cap [-6,0),$$
where $\Ii$ is a DCC set. This is not possible as $(K_X+B+\Mm_X)\cdot C_i$ is strictly decreasing.

We are left to prove the moreover part. Since MMP does not decrease discrepancies, $(X', B',\Mm)$ is lc. By Lemma \ref{lem: qfactorial dlt preserved under mmp}, $X'$ is $\Qq$-factorial klt. Since $\phi$ is $(P+\mu N)$-trivial, $\phi_*(P+\mu N)$ is nef$/U$. We may let $P':=\phi_*(P+\mu N)$ and $N':=(1-\mu)\phi_*N$.
\end{proof}

\begin{prop}\label{prop:log minimal model exist by weak Zariski decomposition}
    Let $(X, B,\Mm)/U$ be an NQC lc generalized pair of dimension $3$ such that $K_X+B+\Mm_X$ admits a weak Zariski decomposition$/U$. Then $(X, B,\Mm)/U$ has a log minimal model.
\end{prop}

\begin{proof}

\noindent\textbf{Step 1}. In this step, we introduce an invariant $\theta(X, B,\Mm,f, N)$. 

Let $\mathfrak{W}$ be the set of NQC lc generalized pairs $(X, B,\Mm)/U$ of dimension $3$ such that $K_X+B+\Mm_X$ admits a weak Zariski decomposition$/U$ but $(X, B,\Mm)/U$ does not have a log minimal model. It suffices to prove $\mathfrak{W}$ is empty. For any $(X, B,\Mm)/U\in\mathfrak{W}$ with weak Zariski decomposition $(f, P, N)$, we let $\theta(X, B,\Mm, f, N)$ be the number of components of $f_*N$ that are not contained in $\lfloor B\rfloor$.

\medskip

\noindent\textbf{Step 2}. In this step we reduce to the case when $\theta(X,B,\Mm,f,N)$ is minimal, $\Supp N\subset\lfloor B\rfloor$, $f=\id_X$, and $(X,B)$ is log smooth, and $\Mm$ descends to $X$.

Assume that $\theta(X,B,\Mm,f,N)$ is minimal, i.e. for any $(X',B',\Mm')/U'\in\mathfrak{W}$ with weak Zariski decomposition $(f',P',N')$, we have $\theta(X',B',\Mm',f',N')\geq\theta(X,B,\Mm,f,N)$. Possibly replacing $f$, we may assume that $f: W\rightarrow X$ is a log resolution of $(X,B+f_*N)$ and $\Mm$ descends to $W$. We let $B_W:=f^{-1}_*B+\Exc(f)$. Then
$$K_W+B_W+\Mm_W=f^*(K_X+B+\Mm_X)+F=P+(N+F)$$
for some $F\geq 0$ that is exceptional$/X$. Then $(\id_W,P,N+F)$ is a weak Zariski decomposition$/U$ of $(W,B_W,\Mm)/U$. Since $F$ is exceptional$/X$ and $\Supp B_W$ contains $\Exc(f)$, we have
$\theta(X,B,\Mm,f,N)\geq\theta(W,B_W,\Mm,\id_W,N+F)$. Since $(X, B,\Mm)/U$ does not have a log minimal model, by Lemma \ref{lem: log minimal model for log smooth model}, $(W, B_W,\Mm)/U$ does not have a log minimal model. Therefore, possibly replacing $(X,B,\Mm)$ with $(W,B_W,\Mm)$, we may assume that $f=\id_X$, $(X,B)$ is log smooth, and $\Mm$ descends to $X$. In particular, $(X,B,\Mm)$ is $\Qq$-factorial dlt.

\medskip

\noindent\textbf{Step 3}. In this step, we deal with the case when $\theta(X, B,\Mm,f, N)=0$. 

Suppose that $\theta(X,B,\Mm,f,N)=0$, then $\Supp N\subset\lfloor B\rfloor$. By Definition-Lemma \ref{deflem: MMP-Zariski}, we may run a $(K_X+B+\Mm_X)$-MMP$/U$ using weak Zariski decomposition $(P,N)$. By Proposition \ref{prop: special termination}, the MMP terminates near a neighborhood of the image of $\Supp \lfloor B \rfloor$. Since $\Supp N\subset\lfloor B\rfloor$, the MMP terminates near the image of $\Supp N$. This is not possible, because for any sequence of steps of  $(K_X+B+\Mm_X)$-MMP$/U$ using weak Zariski decomposition $(P, N)$, then the image of $N$ intersects the extremal ray in the MMP negatively. Thus, the MMP terminates, which contradicts that $X \in \mathfrak{W}$.

\medskip

\noindent\textbf{Step 4}. From now on, we suppose $\theta(X, B,\Mm,f, N)>0$. In this step, we construct an auxiliary generalized pair $(X,\Delta,\Mm)/U$ which admits a weak Zariski decomposition$/U$.

We define 
$$\alpha = \min \{t > 0 | \lfloor (B + t N)^{\leq 1} \rfloor \neq \lfloor B \rfloor \}.$$  
We may write $(B+\alpha N)^{\leq 1}:=B+C$ where $C\geq 0$ and $\Supp C\subset\Supp N$, then $\theta(X,B,\Mm,f,N)$ is the number of components of $C$ and $C\wedge\lfloor B\rfloor=0$. We may write $\alpha N=A+C$, where $A\geq 0$ and $\Supp A\subset\lfloor B\rfloor$. Since $(X,B)$ is log smooth and $\Mm$ descends to $X$, $(X,\Delta:=B+C,\Mm)/U$ is lc, and $(\id_X,P,N+C)$ is a weak Zariski decomposition of $(X,\Delta,\Mm)/U$. By our construction,
$$\theta(X,B,\Mm,f,N)>\theta(X,\Delta,\Mm,f,N+C).$$
Thus by the minimality of $\theta(X,B,\Mm,f,N)$, there exists a log minimal model $(Y,\Delta_Y,\Mm)/U$ of $(X,\Delta,\Mm)$.

Let $g: V \to X$ and $h: V \to Y$ be a common resolution, $P':=h^*(K_Y+\Delta_Y+\Mm_Y)$, and
$$N':=g^*(K_X+\Delta+\Mm_X)-h^*(K_Y+\Delta_Y+\Mm_Y).$$
Then $N'\geq 0$, so $(g,P',N')$ is a weak Zariski decomposition$/U$ of $K_X+\Delta+\Mm_X$. Moreover, since
$$g^*(P+N+C)=P'+N',$$
$g^*(N+C)-N'=P'-g^*P$ is anti-nef$/Y$. Since $N'$ is exceptional$/Y$, $h_*(g^*(N+C)-N')\geq 0$, so by the negativity lemma, $g^*(N+C)-N'\geq 0$. Therefore, $\Supp N'\subset\Supp g^*N$.

\medskip

\noindent\textbf{Step 5}. In this step, we show that $C$ is exceptional$/Y$.

We have
\begin{align*}
    (1+\alpha) g^*(K_X + B + \Mm_X) & = g^*(K_X + B + \Mm_X) + \alpha g^* P + \alpha g^* N \\ & = g^*(K_X + B + \Mm_X) + \alpha g^* P + g^* (A + C) \\ & = (P' + \alpha g^* P) + (N' + g^* A).
\end{align*}
Let $P'':=\frac{1}{1+\alpha}(P' + \alpha g^* P)$ and $N'':=\frac{1}{1+\alpha}(N' + g^* A)$. Then $(g,P'',N'')$ is a weak Zariski decomposition$/U$ of $K_X + B + \Mm_X$. Since $\Supp N'' = \Supp(N' + g^* A) \subseteq \Supp g^* N$, we have $\Supp g_* N'' \subseteq \Supp N$. By the minimality of $\theta(X,B,\Mm,f,N)$, we have $\theta(X,B,\Mm,f,N) = \theta(X,B,\Mm,g,N'')$. Therefore, every component of $C$ is contained in $\Supp g_*N''$. Since $N'$ is exceptional$/Y$, $N''=N'+g^*A$, and no component of $C$ is contained in $\Supp A$, $C$ is exceptional$/Y$. In particular, $(Y,\Delta_Y,\Mm)/U$ is a log birational model of $(X,B,\Mm)/U$.

\medskip

\noindent\textbf{Step 6}. Let $G:=g^*C\wedge N'$, $\Tilde{C}:= g^* C - G$, and $\Tilde{N}':= N' - G$. In this step we show that $\Tilde C$ is not exceptional$/X$.

Assume $\Tilde{C}$ is exceptional over $X$. By our construction,
$$\Tilde{N}' - \Tilde{C}= N'-g^*C'=g^*(K_X + B + \Mm_X) - P'$$ 
is anti-nef$/X$. Since $\Tilde{C}$ is exceptional$/X$, $g_*(N'-g^*C)\geq 0$. By the negativity lemma, $\Tilde N'-\Tilde C\geq 0$, so $\Tilde C=0$. Thus
$$g^*(K_X + B + \Mm_X) - h^*(K_Y +\Delta_Y + \Mm_Y)=\Tilde N'\geq 0,$$
so $(Y,\Delta_Y,\Mm)/U$ is a weak lc model of $(X,B,\Mm)/U$. By Lemma \ref{lem: bs wlcm implies lm}, $(X, B,\Mm)/U$ has a log minimal model, a contradiction.

\medskip

\noindent\textbf{Step 7}. In this step, we conclude the proof. By \textbf{Step 6}, $\Tilde C$ is not exceptional$/X$. Let $\beta > 0$ be the smallest real number such that $\beta N-g_*\Tilde C\geq 0$ and let $\Tilde{A} := \beta g^* N - \Tilde{C}$. Then there exists a component $D$ of $g_*\Tilde{C}$ such that $D\not\subset\Supp g_*\Tilde A$. Thus
\begin{align*}
    (1+\beta) g^*(K_X + B + \Mm_X) & = g^*(K_X + B + \Mm_X) + \beta g^* P + \beta g^* N \\ & = g^*(K_X + B + \Mm_X) + \beta g^* P + \Tilde{A} + \Tilde{C} \\ & = (P' + \beta g^* P) + (\Tilde{N'} + \Tilde{A}).
\end{align*}
Let $P''':=\frac{1}{1+\beta}(P' + \beta g^* P)$ and let $N''':=\frac{1}{1+\beta}(\Tilde{N'} + \Tilde{A})$. Then $g_*N'''\geq 0$ and $N'''$ is anti-nef$/X$, so by the negativity lemma, $N'''\geq 0$. Thus $(g,P''',N''')$ is a weak Zariski decomposition$/U$ of $K_X+B+\Mm_X$ and $\Supp g_*N'''\subset\Supp N$. Since $g_*\Tilde C\leq C$, $D$ is not a component of $g_*\Tilde N'$, so $D$ is not a component of $g_*N'''$. Thus $\theta(X,B,\Mm,f,N)>\theta(X,B,\Mm,g,N''')$, which contradicts the minimality of $\theta(X,B,\Mm,f,N)$.
\end{proof}

\subsection{Existence of log minimal model}

\begin{lem}\label{lem:Bir16 lem 9.2 generalized pair}
    Let $(X, B,\Mm)/U$ be an NQC dlt generalized pair of dimension $\leq 3$ and let $A$ be an ample$/U$ $\Rr$-divisor on $X$.
    Then there exists an effective $A' \sim_{\Rr,U} A$ such that $(X, B+A',\Mm)$ is dlt.
\end{lem}
This lemma is \cite[Lemma 9.2]{Bir16}.
Here, we provide a proof in the setting of generalized pairs for completeness.
\begin{proof}
    Since $(X, B,\Mm)$ is dlt, we can take a log resolution $f: W \to X$ of $(X, B,\Mm)$ such that $\Mm$ descends to $W$ and any $f$-exceptional prime divisor on $W$ has positive log discrepancy. Write $K_W+B_W+\Mm_W = f^*(K_X+B+\Mm_X)$. Take an $f$-anti-ample exceptional divisor $E$ on $W$. For a suitable $0<\varepsilon \ll 1$, we can find $A_W \sim_{\Rr,U} f^*A-\varepsilon E$ such that $\Supp (B_W+A_W+\varepsilon E)$ is snc and the coefficients of the $f$-exceptional divisor in $B_W+A_W+\varepsilon E$ are strictly less than $1$. Let $A' := f_*A_W$. Then we have $f^*(K_X+B+A'+\Mm_X) = K_W+B_W+A_W+\varepsilon E+\Mm_W$. And hence $(X, B+A',\Mm)$ is dlt.
\end{proof}

\begin{lem}\label{thm: existence of minimal model dlt case}
    Let $(X, B,\Mm)/U$ be a $\Qq$-factorial NQC dlt generalized pair of dimension $\leq 3$ and suppose that $K_X+B+\Mm_X$ is pseudo-effective. 
    Then $(X, B,\Mm)/U$ has a log minimal model.
\end{lem}
\begin{proof}
    The proof follows the idea of \cite[Theorem 6.6]{HNT20}.

    \medskip
    
    \noindent\textbf{Step 1}. In this step, we find an effective and big $\Rr$-divisor $H$ and a sequence of $(K_X+B+\Mm_X)$-MMP with scaling of $H$ satisfying (1)-(8).

    By Lemma \ref{lem:Bir16 lem 9.2 generalized pair}, there exists an effective and ample $\Rr$-divisor $H$ on $X$ such that 
    \begin{itemize}
        \item[(1)] $K_X+B+\Mm_X+H$ is nef;
        \item[(2)] $(X,B+H,\Mm)$ is dlt.
    \end{itemize}
    Since $K_X+B+\Mm_X$ is pseudo-effective, we have 
    \begin{itemize}
        \item [(3)] $K_X+B+\Mm_X+H$ is big.
    \end{itemize}
    Running a $(K_X+B+\Mm_X)$-MMP with scaling of $H$, we get a sequence 
    \[ (X,B,\Mm)=(X_0,B_0,\Mm) \dashrightarrow (X_1,B_1,\Mm) \dashrightarrow \cdots (X_n,B_n,\Mm) \dashrightarrow \cdots. \]
    Let 
    \[ \lambda_n = \inf \{t \in \Rr | K_{X_n}+B_n+\Mm_{X_n}+tH_n \text{ is nef} \}\]
    be the scaling coefficients, where $H_n$ is the strict transformation of $H$ on $X_n$.
    If $\lim_{n \to \infty} \lambda_n = \lambda>0$, the above sequence becomes a $(K_X+B+\Mm_X+\lambda H)$-MMP.
    By Lemma \ref{lem: hl22 3.4 char p}, there exists a klt pair $(X,\Delta)$ such that the sequence becomes a $(K_X+\Delta)$-MMP and hence terminates.
    Then we get a minimal model by Lemma \ref{lem: bs wlcm implies lm}.
    Hence, we may assume that
    \begin{itemize}
        \item[(4)] the sequence is infinite;
        \item[(5)] $\lim_{n \to \infty} \lambda_n = 0$.
    \end{itemize}
    Moreover, after truncating and replacing $H$ by $\lambda_0 H_0$, we may assume that. 
    \begin{itemize}
        \item[(6)] $\lambda_0 = 1$ and $\lfloor B+H \rfloor = \lfloor B \rfloor$;
        \item[(7)] the sequence consists of flips;
        \item[(8)] the sequence is isomorphic near $\Supp \lfloor B \rfloor$.
    \end{itemize}

    \medskip
    
    \noindent\textbf{Step 2}. In this step, we construct a generalized dlt pair $(W,G+D_{\infty},\Mm)$ with a birational morphism $\mu:W \to X$, $G = \mu^*\lfloor B \rfloor$ and $(W,D_{\infty},\Mm)$ is $\Qq$-factorial terminal.
    It will satisfy (i)-(iv).

    Note that $(X_n, B_n-\lfloor B_n \rfloor+\lambda_nH_n,\Mm)$ is klt by (2) and (6).
    Let $\mu_n:(W_n,D_n,\Mm) \to (X_n,B_n-\lfloor B_n \rfloor+\lambda_nH_n,\Mm)$ be a terminalization by Lemma \ref{cor: extistence of terminalization}.
    Let $g_{n+1}: W_n\dashrightarrow W_{n+1}$ be the induced birational map.
    Since 
    \[ a(E,X_n,B_n-\lfloor B_n \rfloor+\lambda_nH_n,\Mm) \leq a(E, X_{n+1},B_{n+1}-\lfloor B_{n+1} \rfloor+\lambda_{n+1}H_{n+1},\Mm) \]
    for any exceptional prime divisor $E$ over $X$,
    we see that $g_{n+1}^{-1}:W_{n+1} \dashrightarrow W_n$ does not contract any divisor.
    Note that $X_n \dashrightarrow X_{n+1}$ is small by (7), hence $g_{n+1}$ can only contract the $\mu_n$-exceptional divisors.
    Since there are only finitely many $\mu_0$-exceptional divisors, after truncating, we can assume that 
    \begin{itemize}
        \item[(i)] $g_{n+1}:W_n \dashrightarrow W_{n+1}$ is isomorphic in codimension $1$ for all $n$.
    \end{itemize}
    
    Denote by $h_n: W_0 \dashrightarrow W_n$ the induced birational map.
    For any divisor $E$ on $W=W_0$ which is exceptional over $X$, the coefficient of $E$ in $D_n$ is $1-a(E,X_n,B_n-\lfloor B_n \rfloor+\lambda_nH_n,\Mm)$.
    Then, by the inequality above, we have    
    \[ D_0 \geq (h_1^{-1})_*D_1 \geq (h_2^{-1})_*D_2 \geq \cdots \geq (h_n^{-1})_*D_n \geq \cdots \geq 0\]
    on $W=W_0$.
    Hence the limit $D_{\infty} = \lim_{n \to \infty} (h_n^{-1})_*D_n$ exists.
    Note that we have $D_{\infty} \leq D_0$, hence $(W,D_{\infty},\Mm)$ is also terminal.
    Let $G = \mu_0^*\lfloor B \rfloor$ on $W$.
    We have 
    \begin{itemize}
        \item[(ii)] $K_W+G+D_{\infty}+\Mm_W = \lim_{n\to \infty} (h_n^{-1})_*(\mu^*_n(K_{X_n}+B_n+\lambda_nH_n+\Mm_{X_n}))$.
    \end{itemize}

    Let $\psi_n:Y\to W$ and $\varphi_n: Y \to W_n$ be a common resolution of $h_n: W \dashrightarrow W_n$.
    Note that $K_{X_n}+B_n+\Mm_{X_n}+\lambda_nH_n$ is nef, and hence so is $\mu^*_n(K_{X_n}+B_n+\lambda_nH_n+\Mm_{X_n})$.
    Then applying the negativity lemma for $\mu_0\circ\psi_n$ yields
    \[ \psi_n^*(\mu_0^*(K_X+B+\lambda_nH+\Mm_X)) \geq \varphi_n^*(\mu^*_n(K_{X_n}+B_n+\lambda_nH_n+\Mm_{X_n}))\]
    on $Y$.
    Pushing forward to $W$, we get $\mu_0^*(K_X+B+\lambda_nH+\Mm_X) \geq (h_n^{-1})_*(\mu^*_n(K_{X_n}+B_n+\lambda_nH_n+\Mm_{X_n}))$.
    Note that $\lambda_n \to 0$ by (5), taking the limit, we have 
    \begin{itemize}
        \item[(iii)] $\mu_0^*(K_X+B+\Mm_X) \geq K_W+G+D_{\infty}+\Mm_W$;
    \end{itemize}
    and hence
    \begin{itemize}
        \item[(iv)] $(W,G+D_{\infty},\Mm)$ is dlt.
    \end{itemize}
    
    \medskip
    
    \noindent\textbf{Step 3}. In this step, we show that any $(K_W+G+D_{\infty}+\Mm_W)$-MMP will be isomorphic near $\Supp G$.

    Let $\psi_n:Y\to W$ and $\varphi_n: Y \to W_n$ be a common resolution of $h_n: W \dashrightarrow W_n$.
    Then the negativity lemma yields
    \[ \psi_n^*((h_n^{-1})_*(\mu^*_n(K_{X_n}+B_n+\lambda_nH_n+\Mm_{X_n}))) = \varphi_n^*(\mu^*_n(K_{X_n}+B_n+\lambda_nH_n+\Mm_{X_n})) + F_n \]
    for some effective $\Rr$-divisor $F_n$ exceptional over $W$.
    Note that $K_{X_n}+B_n+\lambda_nH_n+\Mm_{X_n}$ is nef and big, and hence it is semi-ample by \cite[Theorem 1.1]{Wal18}.
    Therefore, the stable base locus of $(h_n^{-1})_*(\mu^*_n(K_{X_n}+B_n+\lambda_nH_n+\Mm_{X_n}))$ is equal to $\psi_n(\Supp (F_n))$.
    
    On the other hand, since the MMP in \textbf{Step 1} is isomorphic near $\Supp \lfloor B \rfloor$, 
    there exists an open subset $U_n \subset X_n$ containing $\Supp \lfloor B \rfloor$ and $U_0 \to U_n$ is an isomorphism.
    Restricting the equation to $\psi_n^{-1}(\mu_0^{-1}(U_0))=\varphi_n^{-1}(\mu_n^{-1}(U_n))$, again by negativity lemma,
    we see that $F_n|_{\psi_n^{-1}(\mu_0^{-1}(U_0))} = 0$.
    Hence, we have $\mu^{-1}_0(U_0) \cap \psi_n(\Supp F_n) = \emptyset$.
    In particular, $\Supp G \cap \psi_n(\Supp F_n) = \emptyset$.

    Note that by (ii), the stable base locus of $K_W+G+D_{\infty}+\Mm_W$ is contained in $\bigcup_{n}\psi_n(\Supp F_n) $, and hence is disjoint with $\Supp G$.
    Then any $(K_W+G+D_{\infty}+\Mm_W)$-MMP will be isomorphic near $\Supp G$.
        
    \medskip
    
    \noindent\textbf{Step 4}. Finally, we conclude the proof.

    By (iv), we can run a $(K_W+G+D_{\infty}+\Mm_W)$-MMP. \textbf{Step 3} yields this MMP is isomorphic near $\Supp G$.
    Hence it is also a $(K_W+D_{\infty}+\Mm_W)$-MMP.
    Since $(W,D_{\infty},\Mm)$ is $\Qq$-factorial terminal, by Theorem \ref{thm: termination of terminal case}, the MMP terminates.
    Hence by Lemma \ref{lem: bir12 2.6}, we have a weak Zariski decomposition of $K_W+G+D_{\infty}+\Mm_W$.
    By (iii), this gives a weak Zariski decomposition of $\mu_0^*(K_X+B+\Mm_X)$.
    Then the conclusion follows from Proposition \ref{prop:log minimal model exist by weak Zariski decomposition}.  
\end{proof}

\begin{thm}\label{thm: existence of minimal model lc case}
    Let $(X, B,\Mm)/U$ be an NQC lc generalized pair of dimension $\leq 3$ and $K_X+B+\Mm_X$ is pseudo-effective. 
    Then $(X, B,\Mm)$ has a log minimal model.
\end{thm}
\begin{proof}
    By Definition-Lemma \ref{deflem: dlt model}, there exists a dlt modification $(Y,B_Y,\Mm) \to (X,B,\Mm)$ over $U$.
    Since $(Y, B_Y,\Mm)$ has a log minimal model, $K_Y+B_Y+\Mm_Y$ has a weak Zariski decomposition.
    Hence $K_X+B+\Mm_X$ has a weak Zariski decomposition.
    Therefore, by Proposition \ref{prop:log minimal model exist by weak Zariski decomposition}, $(X, B,\Mm)/U$ has a log minimal model.
\end{proof}

\section{Termination of pseudo-effective flips}

The goal of this section is to prove the termination of flips for generalized pairs when the generalized log canonical $\Rr$-divisor $K_X+B+\Mm_X$ is relatively pseudo-effective:

\begin{thm}[$=$Theorem \ref{thm: tof intro}]\label{thm: termination of flips general}
    Let $(X, B,\Mm)/U$ be an NQC lc generalized pair of dimension $3$ such that $K_X+B+\Mm_X$ is pseudo-effective$/U$. Then any sequence of $(K_X+B+\Mm_X)$-flips$/U$ terminates.
\end{thm}

To prove Theorem \ref{thm: termination of flips general}, we need to prove the following weaker version first:

\begin{thm}\label{thm: termination of flips x klt}
    Let $(X, B,\Mm)/U$ be a $\Qq$-factorial NQC lc generalized pair of dimension $3$ such that $X$ is klt and that $K_X+B+\Mm_X$ is pseudo-effective$/U$. Then any sequence of $(K_X+B+\Mm_X)$-flips$/U$ terminates.
\end{thm}

\begin{lem}\label{lem:dlt modification for MMP}
    Let $(X_1,B_1,\Mm)$ be an NQC lc generalized pair of dimension $3$ and let 
    \[ \xymatrix{
        (X_1,B_1,\Mm) \ar@{-->}[rr]^{\pi_1} \ar[rd]_{\theta_1} && (X_2,B_2,\Mm) \ar@{-->}[rr]^{\pi_2} \ar[rd]_{\theta_2} \ar[ld]^{\theta_1^+} && (X_3,B_3,\Mm) \ar@{-->}[r]^-{\pi_3} \ar[ld]^{\theta_2^+} & \cdots  \\
        & Z_1 & & Z_2 & &
    }\]
    be a sequence of flips$/U$. Assume that either $X_1$ is $\Qq$-factorial klt, or Theorem \ref{thm: termination of flips x klt} holds. Then there exists a commutative diagram 
    \[ \xymatrix{
        (Y_1,B_{Y_1},\Mm) \ar@{-->}[rr]^{\rho_1} \ar[d]^{h_1} && (Y_2,B_{Y_2},\Mm) \ar@{-->}[rr]^{\rho_2} \ar[d]^{h_2} && (Y_3,B_{Y_3},\Mm) \ar@{-->}[r]^-{\rho_3} \ar[d]^{h_3} & \cdots  \\
        (X_1,B_1,\Mm) \ar@{-->}[rr]^{\pi_1} \ar[rd]_{\theta_1} && (X_2,B_2,\Mm) \ar@{-->}[rr]^{\pi_2} \ar[rd]_{\theta_2} \ar[ld]^{\theta_1^+} && (X_3,B_3,\Mm) \ar@{-->}[r]^-{\pi_3} \ar[ld]^{\theta_2^+} & \cdots  \\
        & Z_1 & & Z_2 & &
    }\]
    where, for each $i\geq 1$, the map $\rho_i:Y_i\bir Y_{i+1}$ is a $(K_{Y_i}+B_{Y_i}+\Mm_{Y_i})$-MMP$/Z_i$ and $(Y_i,B_{Y_i},\Mm)$ is a dlt model of $(X_i,B_i,\Mm)$.
    In particular, the induced sequence
    $$(Y_1,B_{Y_1},\Mm)\dashrightarrow (Y_2,B_{Y_2},\Mm)\dashrightarrow (Y_3,B_{Y_3},\Mm)\dashrightarrow\dots$$
is a sequence of steps of a $(K_{Y_1}+B_{Y_1}+\Mm_{Y_1})$-MMP$/U$.
\end{lem}
\begin{proof}
By Definition-Lemma \ref{deflem: dlt model}, there exists a $\Qq$-factorial dlt modification $h_1:(Y_1,B_{Y_1},\Mm) \to (X_1,B_1,\Mm)$.
Suppose that we have already constructed $\Qq$-factorial dlt modifications $h_i: (Y_i,B_{Y_i},\Mm)\rightarrow (X_i,B_i,\Mm)$ for $i\leq n$ and $\rho_i$ for any $i\leq n-1$ which satisfy our requirements. It suffices to construct $h_{n+1}: (Y_{n+1},B_{Y_{n+1}},\Mm)\rightarrow (X_{n+1},B_{n+1},\Mm)$ and $\rho_n$. 

Since $\theta_n$ is a $(K_{X_n}+B_n+\Mm_{X_n})$-flipping contraction$/U$, $K_{X_n}+B_n+\Mm_{X_n}$ is anti-ample$/Z_n$, so there exists an ample$/Z_n$ $\Rr$-divisor $H_n$ on $X_n$ such that $K_{X_n}+B_n+H_n+\Mm_{X_n}\sim_{\mathbb R,Z_n}0$ and $(X_n,B_n+H_n,\Mm)$ is lc. Let $H_{Y_n}:=h_n^*H_n$. 

If Theorem \ref{thm: termination of flips x klt} holds, then we run a $(K_{Y_n}+B_{Y_n}+\Mm_{Y_n})$-MMP$/Z_n$ which terminates with a log minimal model $(Y_{n+1},B_{Y_{n+1}},\Mm)/Z_n$ of $(Y_n,B_{Y_n},\Mm)/Z_n$ with induced birational map $Y_n\dashrightarrow Y_{n+1}$. Since $X_n\dashrightarrow X_{n+1}$ is the ample model$/Z_n$ of $K_{X_n}+B_n+\Mm_{X_n}$, there exists an induced birational morphism $h_{n+1}:Y_{n+1}\rightarrow X_{n+1}$. By the negativity lemma, $K_{Y_{n+1}}+B_{Y_{n+1}}+\Mm_{Y_{n+1}}=h_{n+1}^*(K_{X_{n+1}}+B_{n+1}+\Mm_{X_{n+1}})$, so $h_{n+1}: (Y_{n+1},B_{Y_{n+1}},\Mm)\rightarrow (X_{n+1},B_{n+1},\Mm)$ is a $\Qq$-factorial dlt modification. We may repeat this process, and the lemma follows. In the following, we may assume that $X_1$ is $\Qq$-factorial klt. By Lemma \ref{lem: qfactorial dlt preserved under mmp}, $X_n$ is $\Qq$-factorial dlt.

Let $0<t\ll 1$ be a real number such that $\theta_n$ is $(K_{X_n}+(1-t)(B_n+\Mm_{X_n}))$-negative and $B_{Y_n}'\geq 0$, where
$$K_{Y_n}+B_{Y_n}'+(1-t)\Mm_{Y_n}:=h_n^*(K_{X_n}+(1-t)(B_n+\Mm_{X_n})).$$
Then
$$K_{X_n}+(1-t)(B_n+\Mm_{X_n})\equiv_{Z_n}s(K_{X_n}+B_n+\Mm_{X_n})$$ for some $s>0$, so
$$K_{Y_n}+B_{Y_n}'+(1-t)\Mm_{Y_n}+\frac{s}{2}H_{Y_n}\equiv_{Z_n}\frac{s}{2}(K_{Y_n}+B_{Y_n}+\Mm_{Y_n}).$$
Since $X_n$ is klt, $(X_n,(1-t)B_n,(1-t)\Mm)$ is klt, so $(Y_n,B_{Y_n}',(1-t)\Mm)$ is klt. Since $H_{Y_n}$ is big$/Z_n$ and nef$/Z_n$, there exist $E_n\geq 0$ and ample$/Z_n$ $\Rr$-divisors $A_{n,m}$ on $Y_n$, such that $\frac{s}{2}H_{Y_n}=A_{n,m}+\frac{1}{m}E_n$. Thus for $m\gg 0$, $(Y_n,B_{Y_n}'+\frac{1}{m}E_n,(1-t)\Mm)$ is klt. Fix $m\gg 0$, then by Lemma \ref{lem: hl22 3.4 char p}, there exists a klt pair $(Y_n,\Delta_n)$ such that $$\Delta_n\sim_{\mathbb R,Z_n}B_{Y_n}'+\frac{1}{m}E_n+(1-t)\Mm_{Y_n}+A_{n,m}=B_{Y_n}'+(1-t)\Mm_{Y_n}+\frac{s}{2}H_{Y_n}.$$
By \cite[Theorem 6.11]{HNT20}, we may run a $(K_{Y_n}+\Delta_{n})$-MMP$/Z_n$ with scaling of $\frac{s}{2}H_{Y_n}$ which terminates. This MMP is also a $(K_{Y_n}+B_{Y_n}+\Mm_{Y_n})$-MMP$/Z_n$ with scaling of $H_{Y_n}$, and we let $(Y_{n+1},B_{Y_{n+1}},\Mm)/Z_n$ be the output of this MMP with induced birational map $Y_n\dashrightarrow Y_{n+1}$. Since $X_n\dashrightarrow X_{n+1}$ is the ample model$/Z_n$ of $K_{X_n}+B_n+\Mm_{X_n}$, there exists an induced birational morphism $h_{n+1}:Y_{n+1}\rightarrow X_{n+1}$. By the negativity lemma, $K_{Y_{n+1}}+B_{Y_{n+1}}+\Mm_{Y_{n+1}}=h_n^*(K_{X_n}+B_n+\Mm_{X_n})$, so $h_{n+1}: (Y_{n+1},B_{Y_{n+1}},\Mm)\rightarrow (X_{n+1},B_{n+1},\Mm)$ is a $\Qq$-factorial dlt modification. We may repeat this process, and the lemma follows. 
\end{proof}

\begin{thm}\label{thm:special termination for lc g-pair}
    Let $(X, B,\Mm)/U$ be an NQC lc generalized pair of dimension $3$. Assume that either $X_1$ is $\Qq$-factorial klt, or Theorem \ref{thm: termination of flips x klt} holds. Then for any sequence of $(K_X+B+\Mm_X)$-flips$/U$
    $$(X,B,\Mm)=:(X_1,B_1,\Mm)\dashrightarrow (X_2,B_2,\Mm)\dashrightarrow\dots\dashrightarrow (X_i,B_i,\Mm)\dashrightarrow\dots,$$
the flipping locus does not intersect $\Nklt(X_i,B_i,\Mm)$ for any $i\gg 0$.
\end{thm}
\begin{proof}
    By Lemma \ref{lem:dlt modification for MMP}, there exists a commutative diagram
    \[\xymatrix{
        (Y_1,B_{Y_1},\Mm) \ar@{-->}[rr]^{\rho_1} \ar[d]^{h_1} && (Y_2,B_{Y_2},\Mm) \ar@{-->}[rr]^{\rho_2} \ar[d]^{h_2} && (Y_3,B_{Y_3},\Mm) \ar@{-->}[r]^-{\rho_3} \ar[d]^{h_3} & \cdots  \\
        (X_1,B_1,\Mm) \ar@{-->}[rr]^{\pi_1} \ar[rd]_{\theta_1} && (X_2,B_2,\Mm) \ar@{-->}[rr]^{\pi_2} \ar[rd]_{\theta_2} \ar[ld]^{\theta_1^+} && (X_3,B_3,\Mm) \ar@{-->}[r]^-{\pi_3} \ar[ld]^{\theta_2^+} & \cdots  \\
        & Z_1 & & Z_2 & &
    }\]
    such that 
    $$(Y_1,B_{Y_1},\Mm)\dashrightarrow (Y_2,B_{Y_2},\Mm)\dashrightarrow\dots\dashrightarrow (Y_i,B_{Y_i},\Mm)\dashrightarrow\dots,$$
    is a sequence of steps of a $(K_{Y_1}+B_{Y_1}+\Mm_{Y_1})$-MMP, each $(Y_i,B_{Y_i},\Mm)$ is $\Qq$-factorial dlt, and each $h_i: (Y_i,B_{Y_i},\Mm)\rightarrow (X_i,B_{i},\Mm)$ is a $\Qq$-factorial dlt modification of $(X_i,B_i,\Mm)$. By Proposition \ref{prop: special termination}, possibly truncating the sequence, we may assume that the $(K_{Y_1}+B_{Y_1}+\Mm_{Y_1})$-MMP is an isomorphism near $\lfloor B_1\rfloor$. Let $D_{X_i}=K_{X_i}+B_i+\Mm_{X_i}$, $D_{Y_i}=K_{Y_i}+B_{Y_i}+\Mm_{Y_i}$, $V_{X_i}=\Nklt(X_i,B_i,\Mm)$, and $V_{Y_i}=\Nklt(Y_i,B_{Y_i},\Mm)$. In the following, we shall show that $\pi_i$ is an isomorphism near $V_{X_i}$ for any $i$, which will conclude the proof.
    
    Suppose that $\pi_i$ is not an isomorphism near $V_{X_i}$. Then there exists  $x\in \Exc(\theta_i) \cap V_{X_i}$ and a curve $\gamma\subset\Exc(\theta_i)$ such that $x\in\gamma$ and $D_{X_i}\cdot\gamma<0$. Then for any $0\leq H\sim_{\mathbb R,Z_i}D_{X_i}$, $H\cdot \gamma<0$, so $x \in \gamma \subset \Supp H$. Thus $x \in \bigcap \Supp H = \Bb(X_i/Z_i,D_{X_i})$. 
    Since $V_{Y_i}=h^{-1}_i(V_{X_i})$ and 
    $$\Bb(Y_i/Z_i,D_{Y_i})=\Bb(Y_i/Z_i,h_i^*D_{X_i})=h_i^{-1}(\Bb(X_i/Z_i,D_{X_i})),$$ 
    we have $V_{Y_i}\cap \Bb(Y_i/Z_i,D_{Y_i})\neq \emptyset$. 

    Let $p: W\rightarrow Y_i$ and $q: W\rightarrow Y_{i+1}$ be a common resolution and let $y_i\in V_{Y_i}\cap\Bb(Y_i/Z_i,D_{Y_i})$ be a closed point. Since $\rho_i$ is an isomorphism near $V_{Y_i}$, $y_{i+1}:=\rho_i(y_i)$ is well-defined, and $p^{-1}(V_i)=q^{-1}(V_{Y_{i+1}})$.
    Then we have $p^*D_{Y_i}=q^*D_{Y_{i+1}}+E$ for some $E\geq 0$ that is exceptional$/Y_{i+1}$. For any $0\leq G_{i+1}\sim_{\mathbb R,Z_i}D_{Y_{i+1}}$, we let
    $$G_i:=p_*(q^*G_{i+1}+E)\sim_{\mathbb R,Z_i}D_{Y_i},$$
    then $p^*G_i=q^*G_{i+1}+E$. Since $y_i\in V_{Y_i}\cap\Supp G_i$,
    $$p^{-1}(y_i)\subset p^{-1}(V_{Y_i})\cap\Supp p^*G_i=q^{-1}(V_{Y_{i+1}})\cap\Supp (q^*G_{i+1}+E).$$
Suppose that $p^{-1}(y_i)\subset\Supp E$. Then there exists a prime divisor $F$ over $X_i$ such that $\Center_{Y_i}F=y_i$ and $\Center_{Y_i}F\subset\Supp E$. Thus 
$$a(F,Y_i,B_{Y_i},\Mm)=a(F,Y_{i+1},B_{Y_{i+1}},\Mm)-\mult_FE<a(F,Y_{i+1},B_{Y_{i+1}},\Mm)$$
which is not possible as $\rho_i$ is an isomorphism near $y_i$. Therefore, $$p^{-1}(y_i)\cap(q^{-1}(V_{Y_{i+1}})\cap\Supp (q^*G_{i+1}))\not=\emptyset,$$
so $y_{i+1}\in V_{Y_{i+1}}\cap\Supp G_{i+1}$. Thus $y_{i+1}\in V_{Y_{i+1}}\cap\Bb(Y_{i+1}/Z_i,D_{Y_{i+1}})$. This is not possible as $D_{Y_{i+1}}$ is semi-ample$/U$ and $\Bb(Y_{i+1}/Z_i,D_{Y_{i+1}})=\emptyset$.
\end{proof}

\begin{lem}\label{lem:  termination of flips inductive}
    Let $\Ii_0\subset [0,+\infty)$ be a finite set such that $1\in\Ii_0$. Let $(X,B,\Mm)/U$ be a $\Qq$-factorial lc generalized pair of dimension $3$ such that $X$ is klt and let
    $$(X,B,\Mm)=:(X_1,B_1,\Mm)\dashrightarrow (X_2,B_2,\Mm)\dashrightarrow\dots (X_i,B_i,\Mm)\dashrightarrow\dots$$
    be an infinite sequence of $(K_X+B+\Mm_X)$-flips$/U$. Let $D$ be an $\Rr$-divisor on $X$ and $\Nn$ a $\bb$-divisor on $X$, such that $X_i\dashrightarrow X_{i+1}$ is also a $(D_i+\Nn_{X_i})$-flip, where $D_i$ is the image of $D$ on $X_i$. 
    
    Assume that $B,D\in\Ii_0$ and $\Nn\in\Nef(U,\Ii_0)$. Then there exists a $\Qq$-factorial lc generalized pair $(Y,B_Y,\Mm)/U$ of dimension $3$ and $\Rr$-divisor $D_Y$ on $Y$ satisfying the following.
    \begin{enumerate}
        \item $Y$ is klt,
        \item There exists an infinite sequence of $(K_Y+B_Y+\Mm_Y)$-flips$/U$
            $$(Y,B_Y,\Mm)=:(Y_1,B_{Y_1},\Mm)\dashrightarrow (Y_2,B_{Y_2},\Mm)\dashrightarrow\dots (Y_i,B_{Y_i},\Mm)\dashrightarrow\dots,$$
        \item $Y_i\dashrightarrow Y_{i+1}$ is also a $(D_{Y_i}+\Nn_{Y_i})$-flip, where $D_{Y_i}$ is the image of $D_Y$ on $Y_i$. 
        \item $B_Y,D_Y\in\Ii_0$.
        \item $\lct(Y,B_Y,\Mm;D_Y,\Nn)>\lct(X,B,\Mm;D,\Nn)$.
    \end{enumerate}
\end{lem}
\begin{proof}
    \noindent\textbf{Step 1}.
    Let $t_i:=\lct(X_i,B_i,\Mm;D_i,\Nn)$ for each $i$.  Then the sequence of $(K_X+B+\Mm_X)$-flips$/U$ is also a sequence of $(K_X+B+t_iD+\Mm_X+t_i\Nn_X)$-flips$/U$ for any $i$. Since $(X_i,B_i+t_iD_i,\Mm+t_i\Nn)$ is lc, $(X_j,B_j+t_iD_j,\Mm+t_i\Nn)$ is lc for any $j\geq i$. Therefore, $t_{j}\geq t_i$ for any $j\geq i$. Thus, possibly truncating the MMP, by Theorem \ref{thm: acc lct}, we may assume that $t:=t_i$ is a constant. By Theorem \ref{thm:special termination for lc g-pair}, possibly truncating more, we may assume that 
    \[ \Exc(\theta_i)\cap \Nklt(X_i,B_i+tD_i,\Mm+t\Nn) = \emptyset\]
    for any $i$.

    \medskip

    \noindent\textbf{Step 2}. In this step we construct $Y,B_{Y},D_{Y}$ and show that they satisfy (1),(4) and (5).
    
    Let $\Pp:=\Mm+t\Nn$.
    By Lemma \ref{lem:dlt modification for MMP}, there exists a commutative diagram 
    \[ \xymatrix{
        (Y_1,\Delta_{Y_1},\Pp) \ar@{-->}[rr] \ar[d]^{h_1} &&  (Y_{N_2},\Delta_{Y_{N_2}},\Pp) \ar@{-->}[r] \ar[d]^{h_2} & \cdots  \\
        (X_{1},B_{1}+tD_{1},\Mm+t\Nn) \ar@{-->}[rr]^{\pi_1} \ar[rd]_{\theta_1} &&  (X_{2},B_{2}+tD_{2},\Mm+t\Nn) \ar@{-->}[r]^-{\pi_2} \ar[ld]^{\theta_1^+} & \cdots  \\
        & Z_1 & &
    }\]
    where, for each $i\geq 1$, the map $Y_{i}\bir Y_{{i+1}}$ is a $(K_{Y_{i}}+\Delta_{Y_{i}}+\Pp_{Y_{N_i}})$-MMP$/Z_i$ and $(Y_{N_i},\Delta_{Y_{N_i}},\Pp)$ is a dlt model of $(X_{i},B_{i}+tD_{i},\Mm+t\Nn)$. Moreover, $N_1=1$. In particular, the sequence $(Y_{N_i},\Delta_{Y_{N_i}},\Pp)$ comes from a $(K_{Y_1}+\Delta_{Y_1}+\Pp_{Y_1})$-MMP$/U$
    \[ \xymatrix{
        (Y_1,\Delta_1,\Pp) \ar@{-->}[rr]^{\rho_1} \ar[rd]_{\mu_1} && (Y_2,\Delta_{Y_2},\Pp) \ar@{-->}[rr]^{\rho_2} \ar[rd]_{\mu_2} \ar[ld]^{\mu_1^+} && (Y_3,\Delta_{Y_3},\Pp) \ar@{-->}[r]^-{\rho_3} \ar[ld]^{\mu_2^+} & \cdots  \\
        & Z'_1 & & Z'_2 & &
    }.\]
    Possibly truncating the MMP, we may assume that $\rho_i$ is a flip for any $i$.

    For each $i\geq 1$, write $B_i=\sum_k b_{i,k} G_{i,k}$ and $D_i = \sum_k d_{i,k}G_{i,k}$, where each $G_{i,k}$ is either an irreducible component of $B_i$ or an irreducible component of $D_i$. Then $B_i+tD_i = \sum_k (b_{i,k}+t d_{i,k})G_{i,k}$ and $b_{i,k}+td_{i,k} \in [0,1]$ by construction.
    Set 
    \begin{align*}
        B_{Y_{N_i}} := \sum_{k:b_{i,k}+td_{i,k}<1} b_{i,k} (h_i)_*^{-1} G_{i,k}, \quad
        D_{Y_{N_i}} := \sum_{k:b_{i,k}+td_{i,k}<1} d_{i,k} (h_i)_*^{-1} G_{i,k}.
    \end{align*}
Then $B_{Y_{N_i}},D_{Y_{N_i}} \in \Gamma_0$ for any $i$. By Definition-Lemma \ref{deflem: dlt model}, 
    \begin{align*}
        \Delta_{Y_{N_i}} &= (h_i)_*^{-1}\left( (B_i+t D_i)^{<1} \right) + (h_i)_*^{-1}\left( (B_i+t D_i)^{=1} \right) + E_i  \\
        &= B_{Y_{N_i}} + t D_{Y_{N_i}} + \lfloor \Delta_{Y_{N_i}} \rfloor,
    \end{align*} 
    where $E_i=\Exc(h_i)$. For different $i$ and $j$, $B_i,D_i$ are strict transforms of $B_j,D_j$ relatively, hence $B_{Y_{N_i}},D_{Y_{N_i}}$ are the strict transforms of $B_{Y_{N_j}},D_{Y_{N_j}}$ relatively. We let $B_{Y_r},D_{Y_r}$ be the birational transforms of $B_{Y_{N_i}},D_{Y_{N_i}}$ on $Y_r$ for any $r$. It is clear that this definition does not depend on the choice of $i$ as $Y_i$ and $Y_j$ are isomorphic in codimension $1$ for any $i,j$.

    For any $r$, $(Y_r,\Delta_{Y_r},\Pp)$ is $\Qq$-factorial dlt, hence $(Y_r,\Delta_{Y_r}^{<1}=B_{Y_r}+tD_{Y_r},\Pp)$ is klt. Thus $Y$ is klt, and 
    \[ \lct(Y_r,B_{Y_r},\Mm;D_{Y_r},\Nn)> t = \lct(X,B,\Mm;D,\Nn). \]

    \medskip

    \noindent\textbf{Step 3}. In this step, we show that 
    \[ \Exc(\mu_r) \cap \Supp \lfloor \Delta_{Y_r} \rfloor= \emptyset \]
    for all $r$.
    
    Consider the commutative diagram
    \[ \xymatrix@C=0.5cm{
        (Y_{N_i},\Delta_{Y_{N_i}},\Pp) \ar@{-->}[r]^-{\rho_{N_i}} \ar[d]^{h_i} & (Y_{N_i+1},\Delta_{Y_{N_i+1}},\Pp) \ar@{-->}[r]^-{\rho_{N_i+1}} & \dots \ar@{-->}[r] & (Y_{N_{i+1}},\Delta_{Y_{N_{i+1}}},\Pp)\ar[d]^{h_{i+1}}  \\
        (X_i,B_i+tD_i,\Mm+t\Nn) \ar@{-->}[rrr]^{\pi_i} \ar[rd]_{\theta_i} &&&  (X_{i+1},B_{i+1}+tD_{i+1},\Mm+t\Nn) \ar[lld]^{\theta_i^+}  \\
        & Z_i & &
    }.\]
    Let $U_{X_i}:=X_i\backslash\Exc(\theta_i), U_{N_i}:=h_i^{-1}(U_{X_i})\subset Y_{N_i}$, and $T_i=\theta_i(U_{X_i})$.
    Since $\Nklt(X_i,B_i+tD_i,\Mm+t\Nn) \subset U_{X_i}$, $\Nklt(Y_{N_i},\Delta_{Y_{N_i}},\Pp) = \Supp \lfloor \Delta_{Y_{N_i}} \rfloor \subset U_{N_i}$. 
    $\pi_i|_{U_{X_i}}:U_{X_i} \to \pi_i(U_{X_i})$ is an isomorphism and hence $(K_{X_i}+B_i+tD_i+\Mm_{X_i}+t\Nn_{X_i})|_{U_{X_i}}$ is trivially semi-ample over $T_i$.
    Then $(K_{Y_{N_i}}+\Delta_{Y_{N_i}}+\Pp_{Y_{N_i}})|_{U_{N_i}} = (h_i|_{U_{X_i}})^*((K_{X_i}+B_i+tD_i+\Mm_{X_i}+t\Nn_{X_i})|_{U_{X_i}})$ is also semi-ample over $T_i$.

    Suppose we have known that $(K_{Y_{r}}+\Delta_{Y_{r}}+\Pp_{Y_{r}})|_{U_{r}}$ is semi-ample over $T_i$ for some open subset $U_r$ such that the induced map $U_{N_i}\dashrightarrow U_r$ is an isomorphism, and $\Supp \lfloor \Delta_{Y_{r}} \rfloor \subset U_{r}$ for $N_i\leq r<N_{i+1}$. Since $\Exc(\mu_{r})$ is covered by $(K_{Y_{r}}+\Delta_{Y_{r}}+\Pp_{Y_{r}})$-negative curves, $U_{r} \cap \Exc(\mu_{r}) = \emptyset$, hence $\rho_{r}$ is an isomorphism over $U_{r}$.
    Thus $\Supp \lfloor \Delta_{Y_{r}} \rfloor\cap \Exc(\mu_{N_i}) = \emptyset$.
    
    Set $U_{r+1} = \rho_r(U_r)$, then $(K_{Y_{r+1}}+\Delta_{Y_{r+1}}+\Pp_{Y_{r+1}})|_{U_{r+1}}$ is semi-ample$/T_i$ and $\Supp \lfloor \Delta_{Y_{r+1}} \rfloor \subset U_{r+1}$. By induction,
       \[ \Exc(\mu_r) \cap \Supp \lfloor \Delta_{Y_r} \rfloor= \emptyset \]
    for all $r$. Therefore, the $(K_{Y_1}+\Delta_{Y_1}+\Mm_{Y_1})$-MMP in \textbf{Step 2} is also a $(K_{Y_1}+ B_{Y_1} + t D_{Y_1}+\Mm_{Y_1}+t\Nn_{Y_1})$-MMP as $\lfloor\Delta_{Y_r}\rfloor\cdot R=0$ for any extremal ray $R$ that is contracted by $\mu_r$.

    \medskip

    \noindent\textbf{Step 4}. In this step we show that $\rho_r:Y_r \bir Y_{r+1}$ is also a $(K_{Y_r}+ B_{Y_r}+\Mm_{Y_r})$-MMP and a $(D_{Y_r}+\Nn_{Y_r})$-flip.

    Suppose $N_i \leq r < N_{i+1}$. 
    By construction of $D_{Y_{N_i}}$, $h_i^*(D_i+\Nn_{X_i}) = D_{Y_{N_i}}+F_{N_i}'+\Nn_{Y_{N_i}}+F_{N_i}''$, where $\mult_{F_{N_i}'}\Delta_{Y_{N_i}}=1$ and $F_{N_i}''$ is $h_i$-exceptional.
    Then we can write 
    \[ h_i^*(D_i+\Nn_{X_i}) = D_{Y_{N_i}}+\Nn_{Y_{N_i}}+F_{N_i},\quad \Supp F_{N_i} \subset \Supp \lfloor \Delta_{Y_{N_i}} \rfloor. \]
    Since $X_i\dashrightarrow X_{i+1}$ is a $(D_i+\Nn_{X_i})$-flip$/Z_i$,
    \[ K_{X_i}+B_i+tD_i+\Mm_{X_i}+t\Nn_{X_i} \equiv_{\Rr,Z_i} (t+\alpha_i)(D_i+\Nn_{X_i}) \]
    for some $\alpha_i>0$.
    Then 
    \[K_{Y_{N_i}}+\Delta_{Y_{N_i}}^{<1}+\lfloor\Delta_{Y_{N_i}}\rfloor+\Pp_{Y_{N_i}} \equiv_{\Rr,Z_i} (t+\alpha_i)(D_{Y_{N_i}}+\Nn_{Y_{N_i}}+F_{N_i})\]
    and then
    \[K_{Y_{r}}+\Delta_{Y_{r}}^{<1}+\lfloor\Delta_{Y_{r}}\rfloor+\Pp_{Y_{r}} \equiv_{\Rr,Z_i} (t+\alpha_i)(D_{Y_{r}}+\Nn_{Y_{r}}+F_{r}),\]
    where $F_r$ is the strict transform of $F_{N_i}$ on $Y_r$.
    In particular, we have 
    \[K_{Y_{r}}+\Delta_{Y_{r}}^{<1}+\lfloor\Delta_{Y_{r}}\rfloor+\Pp_{Y_{r}} \equiv_{\Rr,Z_r'} (t+\alpha_i)(D_{Y_{r}}+\Nn_{Y_{r}}+F_{r}).\]
    By \textbf{Step 3} we know $\Exc(\mu_r) \cap \Supp \lfloor \Delta_{Y_r} \rfloor= \emptyset$ and then $\lfloor \Delta_{Y_r} \rfloor \equiv_{\Rr,Z_r'} F_r \equiv_{\Rr,Z_r'} 0$.
    Hence we get 
    \[ K_{Y_{r}}+B_{Y_r}+tD_{Y_r}+\Mm_{Y_{r}}+t\Nn_{Y_{r}} \equiv_{\Rr,Z_r'} (t+\alpha_i)(D_{Y_{r}}+\Nn_{Y_{r}}).\]
    Then 
    \[ D_{Y_{r}}+\Nn_{Y_{r}} \equiv_{\Rr,Z_r'} \frac{1}{t+\alpha_i}(K_{Y_{r}}+B_{Y_r}+tD_{Y_r}+\Mm_{Y_{r}}+t\Nn_{Y_{r}})\]
    and 
    \[ K_{Y_{r}}+B_{Y_r}+\Mm_{Y_{r}}\equiv_{\Rr,Z_r'} \frac{\alpha_i}{t+\alpha_i}(K_{Y_{r}}+B_{Y_r}+tD_{Y_r}+\Mm_{Y_{r}}+t\Nn_{Y_{r}}).\]
    This proves the assertion after setting $(Y,B_Y,\Mm):=(Y_1,B_{Y_1},\Mm)$ and $D_Y:=D_{Y_1}$.
\end{proof}

\begin{proof}[Proof of Theorem \ref{thm: termination of flips x klt}]
Let
$$(X,B,\Mm)=:(X_1,B_1,\Mm)\dashrightarrow (X_2,B_2,\Mm)\dashrightarrow\dots (X_i,B_i,\Mm)\dashrightarrow\dots$$
be a sequence of $(K_X+B+\Mm_X)$-flips$/U$. 
By Theorem \ref{thm: existence of minimal model lc case}, $(X,B,\Mm)/U$ has a log minimal model $(X_{\min},B_{\min},\Mm)/U$.
Let $p: W\rightarrow X$ and $q: W\rightarrow X_{\min}$ be a common resolution. By Lemma \ref{lem: bir12 2.6},
$$p^*(K_X+B+\Mm_X)=q^*(K_{X_{\min}}+B_{\min}+\Mm_{X_{\min}})+E$$
for some $E\geq 0$. Let
$$\Nn:=\overline{K_{X_{\min}}+B_{\min}+\Mm_{X_{\min}}},$$
then $\Nn$ is a nef$/U$ $\bb$-divisor. By Theorem \ref{thm: shokurov polytope intro}, $\Nn$ is an NQC$/U$ $\bb$-divisor. Let $D:=p_*E$, then
$$D+\Nn_X=K_X+B+\Mm_X.$$
In particular, $D+\Nn_X$ is $\Rr$-Cartier, and $X_i\dashrightarrow X_{i+1}$ is a $(D_i+\Nn_{X_i})$-flip, where $D_i$ is the image of $D$ on $X_i$.

Let $\Ii_0\subset [0,+\infty)$ be a set such that $B,D\in\Ii_0$ and $\Mm\in\Nef(U,\Ii_0)$. By Lemma \ref{lem:  termination of flips inductive}, if this sequence of flips does not terminate, then we may inductively construct a sequence of generalized pairs $(X^i, B^i,\Mm)/U$ and $\Rr$-divisors $D^i$, such that $B^i, D^i\in\Ii_0$, and $$\lct(X^{i+1}, B^{i+1},\Mm, D^{i+1},\Nn)>\lct(X^i, B^i,\Mm, D^{i},\Nn)$$
for each $i$. This contradicts Theorem \ref{thm: acc lct}.
\end{proof}

\begin{proof}[Proof of Theorem \ref{thm: tof intro}]
Suppose that there exists an infinite sequence
$$(X,B,\Mm)=:(X_1,B_1,\Mm)\dashrightarrow (X_2,B_2,\Mm)\dashrightarrow\dots (X_i,B_i,\Mm)\dashrightarrow\dots$$
of $(K_X+B+\Mm_X)$-flips$/U$. By Lemma \ref{lem:dlt modification for MMP}, there exists a commutative diagram 
    \[ \xymatrix{
        (Y_1,B_{Y_1},\Mm) \ar@{-->}[rr]^{\rho_1} \ar[d]^{h_1} && (Y_2,B_{Y_2},\Mm) \ar@{-->}[rr]^{\rho_2} \ar[d]^{h_2} && (Y_3,B_{Y_3},\Mm) \ar@{-->}[r]^-{\rho_3} \ar[d]^{h_3} & \cdots  \\
        (X_1,B_1,\Mm) \ar@{-->}[rr]^{\pi_1} \ar[rd]_{\theta_1} && (X_2,B_2,\Mm) \ar@{-->}[rr]^{\pi_2} \ar[rd]_{\theta_2} \ar[ld]^{\theta_1^+} && (X_3,B_3,\Mm) \ar@{-->}[r]^-{\pi_3} \ar[ld]^{\theta_2^+} & \cdots  \\
        & Z_1 & & Z_2 & &
    }\]
    where, for each $i\geq 1$, the map $\rho_i:Y_i\bir Y_{i+1}$ is a $(K_{Y_i}+B_{Y_i}+\Mm_{Y_i})$-MMP$/Z_i$ and $(Y_i,B_{Y_i},\Mm)$ is a dlt model of $(X_i,B_i,\Mm)$. Thus 
    $$(Y_1,B_{Y_1},\Mm)\dashrightarrow (Y_2,B_{Y_2},\Mm)\dashrightarrow\dots\dashrightarrow (Y_i,B_{Y_i},\Mm)\dashrightarrow\dots$$
is an infinite sequence of $(K_{Y_1}+B_{Y_1}+\Mm_{Y_1})$-MMP$/U$. Since each $Y_i$ is $\Qq$-factorial klt, this contradicts Theorem \ref{thm: termination of flips x klt}.
\end{proof}

\begin{proof}[Proof of Theorem \ref{thm: emm intro}]
Part (2) follows from Theorem \ref{thm: tof intro}. 

For part (1), let $(X',B',\Mm)$ be a dlt model of $(X,B,\Mm)$. By (2) and Lemma \ref{lem: qfactorial dlt preserved under mmp}, $(X', B',\Mm)/U$ has a log minimal model. Thus $K_{X'}+B'+\Mm_{X'}$ admits a weak Zariski decomposition$/U$, so $K_X+B+\Mm_X$ admits a weak Zariski decomposition$/U$. (1) follows from Proposition \ref{prop:log minimal model exist by weak Zariski decomposition}.
\end{proof}


\begin{thebibliography}{99}

\bibitem[Ale94]{Ale94} V. Alexeev, \textit{Boundedness and $K^2$ for log surfaces}, Internat. J. Math. \textbf{5} (1994), 779--810.

\bibitem[AHK07]{AHK07} V.~Alexeev, C. D.~Hacon, and Y.~Kawamata, \textit{Termination of (many) $4$-dimensional log flips}, Invent. Math. \textbf{168} (2007), no. 2, 433--448.


\bibitem[BF25]{BF25} F. Bernasconi and S. Filipazzi, \textit{Rational points on 3-folds with nef anti-canonical class over finite fields}, Proc. Lond. Math. Soc. (3) \textbf{130} (2025), no. 1, Paper No. e70014, 31.

\bibitem[BS23]{BS23} F. Bernasconi and L. Stigant, \textit{Semiampleness for Calabi-Yau surfaces in positive and mixed characteristic}, Nagoya Math. J., \textbf{250} (2023), 365--384.

\bibitem[BMPSTWW23]{BMPSTWW23} B. Bhatt, L. Ma, Z. Patakfalvi, K. Schwede, K. Tucker, J. Waldron, and J. Witaszek, \textit{Globally + regular varieties and the minimal model program for threefolds in mixed characteristic}, Publ. Math. IHÉS \textbf{138} (2023), 69--227.


\bibitem[Bir12]{Bir12} C. Birkar, \textit{Existence of log canonical flips and a special LMMP}, Publ. Math. IHÉS \textbf{115} (2012), 325--368.

\bibitem[Bir16]{Bir16} C. Birkar, \textit{Existence of flips and minimal models for $3$-folds in char $p$}, Ann. Sci. Éc. Norm. Supér. (4) \textbf{49} (2016), no. 1, 169--212.

\bibitem[BCHM10]{BCHM10}
C. Birkar, P. Cascini, C. D. Hacon and J. M\textsuperscript{c}Kernan, \textit{Existence of minimal models for varieties of log general type}, J. Amer. Math. Soc. \textbf{23} (2010), no. 2, 405--468.

\bibitem[BW17]{BW17} C. Birkar and J.Waldron, \textit{Existence of Mori fibre spaces for 3-folds in char $p$}, Advances in Mathematics, \textbf{313} (2017), 62--101.

\bibitem[BZ16]{BZ16} C. Birkar and D.-Q. Zhang, \textit{Effectivity of Iitaka fibrations and pluricanonical systems of polarized pairs}, Publ. Math. IHÉS \textbf{123} (2016), 283--331.


\bibitem[CS21]{CS21} P.~Cascini and C. Spicer, \textit{MMP for co-rank one foliations on threefolds}, Invent. Math. \textbf{225} (2021), no. 2, 603--690.

\bibitem[CS25]{CS25} P. Cascini and C. Spicer, \textit{On the MMP for rank one foliations on threefolds}, Forum Math. Pi. \textbf{13} (2025), Paper No. e20, 38.

\bibitem[CHLX23]{CHLX23} G. Chen, J. Han, J. Liu, and L. Xie, \textit{Minimal model program for algebraically integrable foliations and generalized pairs}, arXiv:2309.15823. 

\bibitem[CT23]{CT23} G.~Chen and N.~Tsakanikas, \textit{On the termination of flips for log canonical generalized pairs}, Acta. Math. Sin. English Ser. \textbf{39} (2023), 967--994. 

\bibitem[CP18]{CP18} V. Cossart and O. Piltant, \textit{Resolution of singularities of threefolds in positive characteristic, I. Reduction to local uniformization on Artin-Schreier and purely inseparable
coverings}, J. Algebra \textbf{320} (2018), 1051--1082.

\bibitem[CP19]{CP19}  V. Cossart and O. Piltant, \textit{Resolution of singularities of threefolds in positive characteristic. II}, J. Algebra \textbf{321} (2009), 1836--1976.

\bibitem[Cut09]{Cut09} S. D. Cutkosky, \textit{Resolution of singularities for 3-folds in positive characteristic}, Amer. J. Math. \textbf{131} (2009), 59--127.

\bibitem[DW22]{DW22} O. Das and J. Waldron, \textit{On the log minimal model program for threefolds over imperfect fields of characteristic}, J. Lond. Math. Soc. (2) \textbf{106} (2022), no. 4, 3895–3937.

\bibitem[FW23]{FW23} S. Filipazzi and J. Waldron, \textit{Connectedness principle for 3-Folds in characteristic $p>5$}, Michigan Math. J. (2023), Advance Publication 1--27.

\bibitem[HL23]{HL23} C. D. Hacon and J. Liu, \textit{Existence of flips for generalized lc pairs}, Camb. J. Math. \textbf{11} (2023), no. 4, 795--828.  

\bibitem[HW22]{HW22} C. D. Hacon and J. Witaszek, \textit{The minimal model program for threefolds in characteristic 5}, Duke Math. J. \textbf{171} (15 August 2022), no. 11, 2193--2231.

\bibitem[HX13]{HX13} C. D. Hacon and C. Xu, \textit{Existence of log canonical closures}, Invent. Math. \textbf{192} (2013), no. 1, 161--195.

\bibitem[HX14]{HX14} C. D. Hacon and C. Xu, \textit{On the three dimensional minimal model program in positive characteristic}, J. Amer. Math. Soc. \textbf{28} (2015), 711--744.

\bibitem[HL22]{HL22} J.~Han and Z.~Li, \textit{Weak Zariski decompositions and log terminal models for generalized polarized pairs}, Math. Z. \textbf{302} (2022), 707--741.

\bibitem[HLS19]{HLS19} J. Han, J. Liu, and V. V. Shokurov, \textit{ACC for minimal log discrepancies of exceptional singularities}, arXiv:1903.04338.


\bibitem[HH20]{HH20}  K. Hashizume and Z. Hu, \textit{On minimal model theory for log abundant lc pairs}, J. Reine Angew. Math. \textbf{767} (2020), 109--159. 

\bibitem[HNT20]{HNT20} K. Hashizume, Y. Nakamura, and H. Tanaka, \textit{Minimal model program for log canonical threefolds in positive characteristic}, Math. Res. Lett. \textbf{27} (2020), no. 4, 1003--1054.

\bibitem[Kol13]{Kol13} J. Koll\'ar, \textit{Singularities of the minimal model program}, Cambridge Tracts in Math. \textbf{200} (2013), Cambridge Univ. Press. With a collaboration of S\'andor Kov\'acs.

\bibitem[KM98]{KM98} J. Koll\'{a}r and S. Mori, \textit{Birational geometry of algebraic varieties}, Cambridge Tracts in Math. \textbf{134} (1998), Cambridge Univ. Press.

\bibitem[PW22]{PW22} Z. Patakfalvi and J. Waldron, \textit{Singularities of general fibers and the LMMP}, American Journal of Mathematics. \textbf{144} (2022) no. 2, 505--540.

\bibitem[Wal17]{Wal17} J. Waldron \textit{Finite generation of the log canonical ring for 3-folds in char $p$}, Mathematical Research Letters. \textbf{24} (2017), no.3 933--946.

\bibitem[Wal18]{Wal18} J. Waldron, \textit{The LMMP for log canonical 3-folds in characteristic $p > 5$}, Nagoya Mathematical Journal. \textbf{230} (2018), 48--71.

\bibitem[Xu24]{Xu24} Z. Xu, \textit{Note on the three-dimensional log canonical abundance in characteristic $>3$}, Nagoya Math. J. (online, 2024), 1--30.
\end{thebibliography}
\end{document}